\documentclass[final,letterpaper]{siamart251216}
\usepackage{amssymb,mathtools}
\usepackage{booktabs}
\usepackage{array,tabularx}
\usepackage{enumitem}
\usepackage{placeins}
\usepackage{graphicx}
\newcommand{\ii}{\mathrm{i}}

\renewcommand{\Im}{\operatorname{Im}}
\newcommand{\R}{\mathbb{R}}
\newcommand{\Z}{\mathbb{Z}}
\newcommand{\vx}{\mathbf{x}}
\newcommand{\vy}{\mathbf{y}}
\newcommand{\vk}{\mathbf{k}}
\newcommand{\vxi}{\boldsymbol{\xi}}
\newcommand{\dxi}{\Delta\xi}
\newcommand{\ximax}{\xi_{\max}}
\newcommand{\BZ}{\mathcal B}
\newcommand{\vm}{\mathbf{m}}
\newcommand{\vn}{\mathbf{n}}
\newcommand{\vell}{\boldsymbol{\ell}}
\newcommand{\Sp}{\mathcal S}
\newcommand{\Sh}{\widehat{\mathcal S}}
\newcommand{\LL}{\mathcal L}
\newcommand{\rhoh}{\widehat\rho}
\newcommand{\Modes}{\mathcal K}
\newsiamremark{remark}{Remark}

\headers{Kernel Design in Particle Mesh Ewald Methods}{S. Jiang, J. Liang, and L. Lu}

\title{Kernel Design in Particle Mesh Ewald Methods: Splitting,
Interpolation, and Upsampling}

\author{Shidong Jiang\thanks{Center for Computational Mathematics, Flatiron
Institute, Simons Foundation, New York, NY 10010
(\email{sjiang@flatironinstitute.org}, \email{jliang@flatironinstitute.org},
\email{llu@flatironinstitute.org}).}
\and Jiuyang Liang\footnotemark[1]\ \thanks{School of Mathematical Sciences,
Shanghai Jiao Tong University, Shanghai 200240, China.}
\and Libin Lu\footnotemark[1]}

\ifpdf
\hypersetup{
  pdftitle={Kernel Design in Particle Mesh Ewald Methods},
  pdfauthor={S. Jiang, J. Liang, and L. Lu}
}
\fi

\begin{document}
\maketitle

\begin{abstract}
For more than thirty years, particle mesh Ewald (PME) methods have been
the dominant way to compute the electrostatic interactions in molecular
dynamics simulations of biomolecules, electrolytes, and materials. They
are implemented in all major simulation packages, where computing these
interactions takes a large share of the run time. A PME method uses two
kernels. A splitting kernel divides the interaction into a local part,
summed over pairs within a cutoff, and a smooth part, computed on a grid.
A spreading and interpolation (SI) kernel spreads the charges onto this
grid and interpolates the potential and the forces from it. The
established methods use Ewald's Gaussian split with B-spline SI kernels,
and a recent method, Ewald summation with prolates, accelerates them
severalfold by using the first prolate spheroidal wave function (PSWF)
of order zero for both kernels. Neither choice meets two needs of
molecular dynamics: an energy that is smooth in the positions of the
particles, which long simulations need in order to conserve it, and a
low-order SI kernel that stays accurate on coarse grids and gains
accuracy as the grid is refined. Both splits leave a jump in the energy
or the force at the cutoff, B-splines need fine grids, and the PSWF
kernel gains little from refinement. We present a framework in which
both kernels are designed to meet these needs. For the inverse powers from the Coulomb to
the dispersion interaction, our splits make the local part smooth at the
cutoff to any prescribed order while keeping the Fourier concentration
of the smooth part close to the optimal one. We show that with analytic
differentiation, used in most PME codes, the root-mean-square force
error of a continuous SI kernel of fixed shape tends to zero under grid
refinement only if every particle, wherever it lies, gives the grid its
exact charge and dipole moment. Our SI kernel for the Coulomb
interaction ensures this by construction, is twice continuously
differentiable, and is designed to reach the requested accuracy at a
given order on the coarsest possible grid. Compared with the B-spline of
the same order, it needs fewer grid points by a factor of up to nearly
five; on coarse grids, it mostly needs a lower order than the PSWF
kernel.
In simulations of liquid water over 100 million time steps, the new
kernels reduce the energy drift of both the established and the PSWF
kernels by about three orders of magnitude, with a smaller force error
than either.
Existing PME codes can adopt them without changing their structure, by
replacing only the two kernels and the Fourier multiplier computed from
them.
\end{abstract}


\section{Introduction}
\label{sec:intro}

Molecular dynamics simulations of proteins, nucleic acids, electrolytes,
and materials require the electrostatic energy and forces of $N$ point
charges $q_1,\ldots,q_N$, with $\sum_iq_i=0$, at
positions $\vx_1,\ldots,\vx_N$ in a periodic box of side $L$:
\begin{equation}
 E=\frac12\sum_{\vn\in\Z^3}\,{\sum_{i,j=1}^N}{}^{\prime}\,
     \frac{q_iq_j}{|\vx_i-\vx_j+L\vn|},\qquad
 \mathbf F_i=-\nabla_{\vx_i}E,
 \label{eq:coulomb-energy}
\end{equation}
where the prime excludes the terms with $i=j$ when $\vn=\mathbf 0$. The
lattice sum converges only conditionally, and direct summation over pairs
costs $O(N^2)$ operations, which a simulation spends again at each of its
time steps, often millions of them.

Most fast algorithms for these sums belong to one of two classes.
Tree-based methods,
such as the fast multipole method~\cite{GreengardRokhlin} and the
dual-space multilevel kernel-splitting method~\cite{DMK}, achieve linear
complexity at fixed accuracy, even for strongly nonuniform
distributions. Methods based on the fast Fourier transform (FFT) use
the splitting that Ewald~\cite{Ewald} introduced in 1921: the Coulomb
kernel is split into a local part, summed directly over nearby pairs, and
a smooth part, whose Fourier series converges rapidly. Which class is
faster depends on the system and the hardware~\cite{Kohnke}; in
molecular dynamics, where the particles fill the box with a nearly
uniform density, the particle mesh Ewald (PME)
methods~\cite{Darden,Essmann,HE}, developed more than thirty years ago,
dominate.

PME methods, implemented in all major molecular dynamics packages,
spread each charge to nearby grid points with a compactly supported
kernel, transform the grid values by the FFT, scale each Fourier
coefficient by a diagonal multiplier, which the molecular dynamics
literature calls the \emph{influence function}~\cite{HE}, transform back, and
interpolate the potential or the force from the grid with the same
kernel. With a cutoff radius that need not grow with $N$ at fixed
density, the grid has $O(N)$ points and the total cost is $O(N\log N)$.
On parallel computers, however, every FFT requires global
communication, and in the benchmarks of~\cite{ESP}, the computation of
the Coulomb interactions takes from a quarter to more than four fifths of
the run time, more for larger systems and higher accuracy. Reaching the
same accuracy on a grid with fewer points can therefore speed up the
whole simulation severalfold~\cite{ESP}.

The accuracy of these methods is controlled by two kernels. The splitting
kernel determines how fast the Fourier transform of the smooth part
decays, and hence how fine the grid must be for a given cutoff. The
spreading and interpolation (SI) kernel determines the aliasing error: a
grid cannot distinguish a Fourier mode from its shifts by multiples of
the sampling frequency, and the SI kernel controls how much of each
shifted copy enters the result. The influence function divides out the
effect of the SI kernel on the retained modes and can be chosen to
minimize the remaining error~\cite{HE}. Most PME codes compute the forces
by analytic differentiation (AD), which interpolates with the gradient of
the SI kernel and needs one inverse FFT instead of the three of spectral
differentiation~\cite{Ball,DH1,ESP}.

The literature on these methods is vast, and we will not attempt a
comprehensive review here. The particle--particle--particle--mesh method
(P3M)~\cite{DH1,DH2,HE} and the smooth particle mesh Ewald method
(SPME)~\cite{Essmann} use Ewald's Gaussian split and B-spline SI kernels;
P3M chooses the influence function that minimizes the mean-square force
error, and the optimal influence functions for AD and for the energy were
derived later~\cite{Ball,EnergyP3M,Stern}. Other methods use
Gaussian~\cite{SE,GSE} or Kaiser--Bessel~\cite{Gao} SI kernels, optimized
interpolation bases~\cite{Wang}, midtown splines~\cite{Midtown}, the
spreading kernels of nonuniform FFTs~\cite{Nestler,Pippig}, or other
splits~\cite{Natoli,USeries}.

Ewald summation with prolates (ESP)~\cite{ESP}, introduced recently,
replaces both kernels of P3M and SPME by the first prolate spheroidal
wave function (PSWF) of order zero, $\psi_0$, the function, among those
supported on an interval, whose Fourier transform is most concentrated in
a given band~\cite{LandauPollak,Slepian}. As a splitting kernel, it needs nearly
half as many grid points in each dimension as the Gaussian at the same
cutoff and high accuracy; as an SI kernel, it was chosen to minimize the number of grid
points coupled to each particle. In LAMMPS and GROMACS, ESP accelerates
the long-range electrostatic calculation by factors of about three to
ten, and complete simulations by factors of about two and a half to five,
depending on the accuracy~\cite{ESP}. The $\psi_0$ split has since been
extended to the inverse powers $1/r^\nu$ that appear in machine-learning
interatomic potentials, such as the dispersion interaction
$1/r^6$~\cite{MLIP}.

Since $\psi_0$ is optimally concentrated, it was argued in~\cite{ESP}
that no significant kernel-based improvement is possible. For the
truncation error of the Coulomb split, this is exactly right: we show
that in large boxes, $\psi_0$ minimizes this error relative to the force
of the smooth part. But molecular dynamics places two requirements on
the kernels that neither the established kernels nor those of ESP meet.
The first is smoothness: where pairs of particles cross the cutoff, and
where the set of grid points coupled to a particle changes, an energy
that is not smooth enough makes the integrator lose its order of
accuracy, and the resulting errors can accumulate into a drift of the
energy~\cite{CottrellTupper,EngleSkeel}. The truncated Gaussian split
and the $\psi_0$ splits make the energy or the force of a pair jump at
the cutoff, and the $\psi_0$ SI kernel jumps at the ends of its support.
In the simulations of liquid water without a thermostat in
\cref{sec:nve}, the energy drifts steadily with both P3M and ESP.

The second requirement concerns upsampling, the refinement of the grid
at a fixed order of the SI kernel. Molecular dynamics codes keep
the order of the SI kernel (the number of grid points coupled to each
particle in each dimension) low, typically four or five, and refine the
grid to reach higher accuracy~\cite{ESP}; on graphics processing units
(GPUs), GROMACS supports only order four~\cite{GromacsGPU}. A low-order SI kernel should therefore be
accurate on coarse grids, as the $\psi_0$ kernel is from order five on,
and become more accurate as the grid is refined, as B-splines do. Neither kernel does both.
B-splines need fine grids at the orders used in practice, and the
$\psi_0$ kernel gains little from refinement. A B-spline gives the grid
the exact charge of a particle wherever the particle lies, but the
$\psi_0$ kernel, whatever its bandwidth, does not (\cref{fig:charge});
with AD, the variation of the assigned charge produces a spurious force
that grows, relative to the true force, as the grid is refined
(\cref{sec:refinement}).

In this paper, we design the two kernels in turn, the split first,
since it depends only on the interaction, the cutoff, and the requested
accuracy. We seek splits whose local part has a prescribed number of
continuous derivatives at the cutoff and whose smooth part keeps the
concentration of its Fourier transform in the band close to the optimum.
The smoothness conditions at the cutoff are linear in the splitting
kernel, and we meet
them with the fewest and most concentrated even PSWFs, starting with
$\psi_0$. Our examples are the $C^2$ splits, which make the energy of
every pair twice continuously differentiable where the pair crosses the
cutoff.

For the SI kernel, which we design for the Coulomb interaction, we start
from the choice that a code makes: it fixes the order and seeks the
coarsest grid on which the requested accuracy is reached. We measure the
accuracy by the root-mean-square (RMS) force error, the standard measure
in molecular dynamics, which for independent, uniformly distributed
particles is a normalized $L^2$ norm of the force error over their
positions. With AD, this error vanishes under upsampling, for a
continuous kernel of fixed shape, only if the kernel gives the grid the
exact charge and dipole moment of a particle at every position; the
twice continuously differentiable kernels that do so are exactly the
convolutions of the linear B-spline with continuous functions. Among
those of a given order that decrease away from their center and are no
rougher than the B-spline of that order, the optimized kernel is the one
that meets the tolerance on the coarsest grid. Since the B-spline is
among them, the optimized kernel never needs a finer grid than the
B-spline. For every order of four or more, the optimized kernel exists,
and we compute these kernels numerically once and tabulate them. SI kernels have been
optimized before for the Gaussian split~\cite{Wang,Midtown}; what is new
here is the optimization over the smooth kernels that conserve charge
and dipole, for the grid that a code chooses.

The new kernels have the following features.
\begin{enumerate}[leftmargin=*]
\item For every inverse power that we treat, the $C^2$ splits remove
the jumps of the established and the $\psi_0$ splits at the cutoff, at
the price of a mildly finer grid: to reach the
truncation error of the $\psi_0$ Coulomb split, the $C^2$ Coulomb split
needs $9\%$ to $14\%$ more grid points in each dimension.
\item On the coarse grids that ESP uses, the optimized kernel needs a
lower order than the $\psi_0$ kernel at the tolerances from $10^{-3}$ to
$10^{-6}$ and the same order at $10^{-2}$, even when the bandwidth of the
$\psi_0$ kernel is tuned for each grid. B-splines need higher orders
still, above ten from $10^{-4}$ on.
\item Under upsampling, the error of the optimized kernel keeps falling,
and the B-spline of the same order needs up to $3.5$ times as many grid
points at the orders four and five, and up to $4.8$ times as many over
the orders four to ten.
\item Together, the two kernels make the energy twice continuously
differentiable in the positions of the particles. Over $10^8$ time steps
of liquid water without a thermostat, the energy then drifts about three
orders of magnitude more slowly than with P3M or ESP. The new kernels
also give the smallest force error of the three.
\item Existing PME codes can use the new kernels without changing their
structure, as the integration of ESP into LAMMPS and GROMACS
showed~\cite{ESP}: the five steps of a PME method and their order stay
the same, and
only the local part of the interaction, the SI kernel, and the influence
function change. The influence function is computed once for each grid,
and for the Coulomb interaction, it is the optimal one.
\end{enumerate}

In short, the kernels of a PME method should be chosen for the error of
the method, under the smoothness and the conservation properties that
molecular dynamics needs. Smoothness costs the split a slightly finer
grid, and conservation brings the SI kernel large gains.

The new kernels thus improve on both the established kernels of P3M and
SPME and those of ESP, in the accuracy reached on a grid of given size
and in the conservation of energy. With order four or five, the
defaults of molecular dynamics codes, and with at most fourfold
upsampling, the optimized kernel meets $10^{-6}$, whereas B-splines meet
no tolerance below $10^{-5}$ and the $\psi_0$ kernel none below $10^{-3}$.
In complete calculations for independent particles and for liquid water
at the tolerance $10^{-4}$, typical in molecular dynamics~\cite{ESP}, the
RMS force error on the coarsest grid of the optimized kernel is up to
eight times as large with the B-spline of the same order as with the
optimized kernel, and up to fifty times as large with the $\psi_0$
kernel.

In \cref{sec:prelim}, we collect the conventions and known results that
we use. \Cref{sec:split-design} states smoothness at the cutoff as linear
conditions on the splitting kernel, shows that the $\psi_0$
split minimizes the truncation error of the Coulomb split, and constructs
the $C^2$ splits.
\Cref{sec:si} derives the RMS force error in terms of the two kernels,
shows why a continuous SI kernel must conserve charge and dipole, constructs the
SI kernels, and explains how to choose the parameters of the method.
\Cref{sec:numerics} tests the new kernels, including their conservation
of energy over $10^8$ time steps, and \cref{sec:conclusions} draws
conclusions.

\section{Analytical preliminaries}
\label{sec:prelim}

This section collects the conventions and known results used below.
The method depends on two functions that
this paper designs: the splitting kernel $\chi$ of \cref{sec:splits},
defined for the inverse powers by a window function $w$, and the
one-dimensional kernel $\phi$ that defines the SI kernel of
\cref{sec:si-kernels}. Every other quantity below, from the two
parts of the split and the SI kernel $\Phi$ to their Fourier transforms
and the error of the method, is determined by these two functions, the
grid, and the cutoff.

\subsection{Fourier conventions and Poisson summation}
\label{sec:fourier}

For a function $f$ on $\R^d$, we use the Fourier transform and its
inverse in the form
\begin{equation}
 \widehat f(\vxi)=\int_{\R^d}f(\vx)e^{\ii\vxi\cdot\vx}\,d\vx,\qquad
 f(\vx)=\frac1{(2\pi)^d}\int_{\R^d}\widehat f(\vxi)
        e^{-\ii\vxi\cdot\vx}\,d\vxi .
 \label{eq:fourier}
\end{equation}
A function on the periodic box $[0,L)^3$ has the Fourier series
\begin{equation}
 f(\vx)=\frac1{L^3}\sum_{\vk\in\Z^3}\widehat f_\vk\,
        e^{-\ii\vxi_\vk\cdot\vx},\qquad
 \widehat f_\vk=\int_{[0,L)^3}f(\vx)e^{\ii\vxi_\vk\cdot\vx}\,d\vx,
 \qquad
 \vxi_\vk=\dxi\,\vk,\quad \dxi=\frac{2\pi}L .
 \label{eq:fourier-series}
\end{equation}

On a grid of spacing $h$, the exponentials $e^{\ii\vxi\cdot\vx}$ and
$e^{\ii(\vxi+2\pi\vm/h)\cdot\vx}$, $\vm\in\Z^d$, take the same values at
the grid points. Poisson summation~\cite{TrapReview} expresses a sum
over the grid points in terms of the frequencies $\vxi+2\pi\vm/h$: for an
even function
$f$ on $\R^d$ with compact support,
\begin{equation}
 \sum_{\vell\in\Z^d}f(\vx-h\vell)\,e^{\ii\vxi\cdot h\vell}
 =\frac1{h^d}\sum_{\vm\in\Z^d}
   \widehat f\Bigl(\vxi+\frac{2\pi\vm}h\Bigr)
   e^{\ii(\vxi+2\pi\vm/h)\cdot\vx}.
 \label{eq:poisson}
\end{equation}
The left side times $e^{-\ii\vxi\cdot\vx}$ is periodic with period $h$
in each coordinate of $\vx$, and the identity follows by expanding this
product in a Fourier series. For a function with jumps, such as the $\psi_0$ kernel
below, the identity holds wherever the left side is continuous. The term with
$\vm=\mathbf 0$ on the right side involves the transform at $\vxi$
itself, and the other terms involve the transform at the aliases
$\vxi+2\pi\vm/h$ of $\vxi$.

\subsection{Prolate spheroidal wave functions}
\label{sec:pswf}

For a bandwidth $c>0$, the PSWFs of order zero, $\psi_n^c$, $n=0,1,2,\ldots$, are the eigenfunctions of the
finite Fourier transform, which maps a function $f$ on $[-1,1]$ to the
function $\int_{-1}^1f(t)e^{\ii cxt}\,dt$ of
$x\in[-1,1]$~\cite{LandauPollak,Osipov,Slepian}. They are real, $\psi_n^c$ has the parity of $n$, and we normalize them so
that
\begin{equation}
 \int_{-1}^1\psi_n^c(t)\psi_m^c(t)\,dt=\delta_{nm},\qquad
 \psi_{2j}^c(0)>0 .
 \label{eq:pswf-normalization}
\end{equation}
Since all kernels in this paper are even, only the even functions
$\psi_{2j}^c$ appear below. Each $\psi_n^c$ is a multiple of its finite
Fourier transform, an entire function of $x$, and therefore extends to an
entire function. The PSWFs are also the eigenfunctions of
the prolate differential equation~\cite{Osipov,Slepian}
\begin{equation}
 \bigl((1-t^2)\psi'(t)\bigr)'+(\lambda_n-c^2t^2)\psi(t)=0,\qquad
 \psi=\psi_n^c,\quad |t|<1,
 \label{eq:prolate-ode}
\end{equation}
whose eigenvalues satisfy $\lambda_0<\lambda_1<\lambda_2<\cdots$, and
$\psi_n^c$ has exactly $n$ zeros in $(-1,1)$~\cite{Osipov}. Since
$\lambda_0$ is the minimum of
$\int_{-1}^1\bigl((1-t^2)\psi'(t)^2+c^2t^2\psi(t)^2\bigr)\,dt
/\int_{-1}^1\psi(t)^2\,dt$ over nonzero $\psi$, and the constant function
gives $c^2/3$, we have $\lambda_0\le c^2/3$. The function $\psi_0^c$ has no zeros in
$[-1,1]$; in particular, $\psi_0^c(1)\ne0$ for every $c>0$~\cite{Osipov}.

Extend $\psi_n^c$ by zero outside $[-1,1]$. Its Fourier transform at
$\omega$ is its finite Fourier transform at $x=\omega/c$, and hence a
multiple of $\psi_n^c(\omega/c)$ for $|\omega|\le c$. Together with
Parseval's identity and \cref{eq:pswf-normalization}, this gives
\begin{equation}
 \int_{-c}^{c}\widehat{\psi_n^c}(\omega)\,
   \overline{\widehat{\psi_m^c}(\omega)}\,d\omega=0\quad(n\ne m),
 \qquad
 \int_{\R}\widehat{\psi_n^c}(\omega)\,
   \overline{\widehat{\psi_m^c}(\omega)}\,d\omega=2\pi\delta_{nm}.
 \label{eq:pswf-double-orthogonality}
\end{equation}
Thus the Fourier transforms of the PSWFs are orthogonal both on the band
$[-c,c]$ and on the whole line. The fraction
of the energy of $\widehat{\psi_n^c}$ in the band decreases with $n$: it
is exponentially close to one for $n$ well below $2c/\pi$ and
exponentially small for $n$ well above~\cite{LandauPollak,Osipov}. Among
all functions supported on $[-1,1]$, the multiples of $\psi_0^c$ alone
have the largest fraction of their energy in the band~\cite{Slepian}. Moreover, for $f=\sum_na_n\psi_n^c$,
\cref{eq:pswf-double-orthogonality} gives
\begin{equation}
 \int_{|\omega|>c}|\widehat f(\omega)|^2\,d\omega
 =\sum_n|a_n|^2\int_{|\omega|>c}|\widehat{\psi_n^c}(\omega)|^2\,d\omega,
 \label{eq:pswf-out-of-band}
\end{equation}
so the energy of a PSWF expansion outside the band is a diagonal
quadratic form in its coefficients, with weights that are exponentially
small for $n$ well below $2c/\pi$.

\subsection{Splitting kernels}
\label{sec:splits}

As in ESP~\cite{ESP}, we split the Coulomb kernel as
$1/r=\Sp(r)+\LL(r)$, where the local part $\LL$ vanishes beyond a cutoff
$r_c$ and the smooth part $\Sp$ has a Fourier transform that is small beyond a
band of frequencies.
The split is determined by a splitting kernel.
For the Coulomb kernel, this is an even function $\chi$ on $[-1,1]$ with
$\int_{-1}^1\chi(t)\,dt=1$, and
\begin{equation}
 \Sp(r)=\frac2r\int_0^{\min(r/r_c,1)}\chi(t)\,dt,\qquad
 \LL(r)=
 \begin{cases}
  \displaystyle\frac2r\int_{r/r_c}^1\chi(t)\,dt, & 0<r<r_c,\\[1.5ex]
  0, & r\ge r_c .
 \end{cases}
 \label{eq:coulomb-split}
\end{equation}
For $\vxi\ne\mathbf0$, the Fourier transform of the smooth part is
\begin{equation}
 \Sh(\vxi)=\frac{4\pi}{|\vxi|^2}\,\widehat\chi(|\vxi|r_c),\qquad
 \widehat\chi(\omega)=\int_{-1}^1\chi(t)e^{\ii\omega t}\,dt .
 \label{eq:split-transform}
\end{equation}
To see this, put $\xi=|\vxi|$. The Fourier transform of a radial
function $v(|\vx|)$ is $(4\pi/\xi)\int_0^\infty rv(r)\sin(\xi r)\,dr$, so
exchanging the order of integration in the transform of $\LL$ gives
$\widehat\LL(\vxi)=(4\pi/\xi^2)(1-\widehat\chi(\xi r_c))$, and
\cref{eq:split-transform} follows by subtracting $\widehat\LL$ from the
transform $4\pi/\xi^2$ of $1/r$. ESP uses the $\psi_0$ split
\begin{equation}
 \chi(t)=\frac{\psi_0^{c_s}(t)}{\int_{-1}^1\psi_0^{c_s}(\tau)\,d\tau},
 \qquad |t|\le1,
 \label{eq:psi0-split}
\end{equation}
where $c_s$ is the splitting bandwidth. Since the finite Fourier
transform of $\psi_0^{c_s}$ is a multiple of $\psi_0^{c_s}$ and
$\widehat\chi(0)=1$, we have $\widehat\chi(c_sx)=\psi_0^{c_s}(x)/\psi_0^{c_s}(0)$ for $|x|\le1$.
Beyond $|\vxi|=c_s/r_c$, the transform \cref{eq:split-transform} is
small because of the concentration of $\psi_0^{c_s}$, and ESP retains the modes of a grid whose Nyquist frequency is close to
$c_s/r_c$. ESP's rule for the splitting bandwidth is
$\chi(1)=\varepsilon$, the requested tolerance for the force~\cite{ESP}.

With $\sum_jq_j=0$ and the conducting boundary condition, the potential
at particle $i$ is
\begin{equation}
 \begin{gathered}
 u_i=\sum_{j=1}^N\,{\sum_{\vn\in\Z^3}}{}^{\prime}\,
      q_j\LL(|\vx_i-\vx_j+L\vn|)+u^s(\vx_i)-q_i\Sp(0),\\
 u^s(\vx)=\frac1{L^3}\sum_{\vk\ne\mathbf0}\Sh(\vxi_\vk)\rhoh(\vxi_\vk)
          e^{-\ii\vxi_\vk\cdot\vx},
 \end{gathered}
 \label{eq:split-potential}
\end{equation}
where the prime excludes $(j,\vn)=(i,\mathbf0)$,
$\rhoh(\vxi)=\sum_jq_je^{\ii\vxi\cdot\vx_j}$, and the last term, with
$\Sp(0)=2\chi(0)/r_c$, removes the interaction of particle $i$ with its
own smooth part. Then $E=\frac12\sum_iq_iu_i$, and the force on particle
$i$ is the sum of the local forces and $-q_i\nabla u^s(\vx_i)$. The local
part is summed directly; the particle--mesh method approximates $u^s$ and
its gradient.

For the inverse powers $r^{-\nu}$, $\nu=1,2,\ldots$, the construction
of~\cite{MLIP} works with the dimensionless radius $t=r/r_c$ and one
window function for all $\nu$: an even function $w$ with $w(0)=1$. The
smooth parts of $t^{-1}$ and $t^{-2}$ are
\begin{equation}
 s_1(t)=\frac{\int_0^tw(\tau)\,d\tau}{t\int_0^1w(\tau)\,d\tau},\qquad
 s_2(t)=\frac{1-w(t)}{t^2},
 \label{eq:inverse-base}
\end{equation}
and those of the higher powers follow from the identity
$-(\nu t)^{-1}\,d(t^{-\nu})/dt=t^{-\nu-2}$:
\begin{equation}
 s_{\nu+2}(t)=-\frac1{\nu t}\,s_\nu'(t),
 \label{eq:inverse-recurrence}
\end{equation}
applied separately on $0<t<1$ and $t>1$. The $\psi_0$ split
of~\cite{MLIP} uses the window function $w(t)=\psi_0^{c_s}(t)/\psi_0^{c_s}(0)$
for every $\nu$, takes its Fourier transform to vanish outside
$[-c_s,c_s]$, so that the smooth parts are band-limited, and truncates the
local parts at the cutoff. We extend $w$ by zero outside $[-1,1]$ instead,
as ESP does for the Coulomb split. Then $s_1$ and $s_2$ equal $t^{-1}$ and
$t^{-2}$, respectively, for $t>1$, and the recurrence gives $s_\nu(t)=t^{-\nu}$ there.
The smooth and local parts of $r^{-\nu}$ are
$\Sp_\nu(r)=r_c^{-\nu}s_\nu(r/r_c)$ and $\LL_\nu(r)=r^{-\nu}-\Sp_\nu(r)$,
so the local parts vanish beyond the cutoff, and the only error of the
split is the truncation error of the Fourier series of its smooth part.
For $\nu=1$, this is the Coulomb split \cref{eq:coulomb-split} with
$\chi=w/\int_{-1}^1w(t)\,dt$, and the $\psi_0$ window gives the split
\cref{eq:psi0-split}. We treat $\nu=1,\ldots,6$, which include the
Coulomb interaction, $\nu=1$, and the dispersion interaction, $\nu=6$.

\subsection{SI kernels}
\label{sec:si-kernels}

Fix a grid spacing $h$ and an integer $P\ge2$, which we call the order of
the SI kernel, as is standard for B-splines~\cite{DH1,Essmann}. The SI
kernel is the tensor product
\begin{equation}
 \Phi(\vx)=\frac1{h^3}\prod_{\alpha=1}^3\phi\Bigl(\frac{x_\alpha}h\Bigr),
 \qquad
 \widehat\Phi(\vxi)=\prod_{\alpha=1}^3\widehat\phi(h\xi_\alpha),\qquad
 \widehat\phi(q)=\int_{\R}\phi(t)e^{\ii qt}\,dt,
 \label{eq:si-kernel}
\end{equation}
where the one-dimensional kernel $\phi$ is even, supported on $[-P/2,P/2]$, and normalized by
$\int_{\R}\phi(t)\,dt=1$, so that $\widehat\Phi(\mathbf 0)=1$. A particle at $\vy$
is thus coupled to the $P^3$ grid points in the support of
$\Phi(\cdot-\vy)$, its stencil. We say that the particle crosses the
boundary of a stencil when, as it moves, a grid point enters or leaves
this support. Kernels that are not tensor products, such as the
cubature splines of~\cite{Midtown}, have stencils of other shapes and
can reach the order of accuracy of the B-spline of order $P$ with fewer
than $P^3$ grid points; we consider only tensor products, which most
molecular dynamics packages use. The two kernels we compare with are the
following.
\begin{itemize}[leftmargin=*]
\item The B-spline $B_P$, the centered cardinal B-spline of order $P$. It
is a piecewise polynomial of degree $P-1$ with $P-2$ continuous
derivatives, and $\widehat\phi(q)=\bigl(\sin(q/2)/(q/2)\bigr)^P$.
\item The $\psi_0$ kernel of ESP,
\begin{equation}
 \phi(t)=\frac{\psi_0^{c_w}(2t/P)}
           {\int_{-P/2}^{P/2}\psi_0^{c_w}(2\tau/P)\,d\tau},\qquad |t|\le P/2,
 \label{eq:psi0-kernel}
\end{equation}
and $\phi(t)=0$ for $|t|>P/2$, where $c_w$ is the SI bandwidth. Since
$\psi_0^{c_w}(1)\ne0$, the $\psi_0$ kernel jumps at the ends of its
support.
\end{itemize}

\begin{figure}[htbp]
\centering
\includegraphics[width=\linewidth]{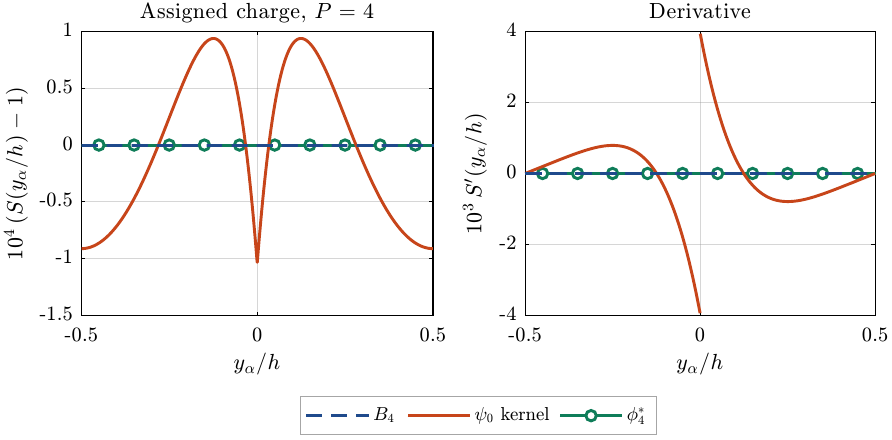}
\caption{The deviation of the factor $S(y_\alpha/h)$ of the assigned
charge \cref{eq:assigned-charge} from one (left), and its derivative
$S'(y_\alpha/h)$ (right), as the particle moves over one grid cell, for
$P=4$: the cubic B-spline $B_4$, the $\psi_0$ kernel for the tolerance
$10^{-4}$ on the grid that ESP uses (\cref{sec:band-results}), and the
optimized kernel $\phi_4^*$ of \cref{sec:si-design} for the same
tolerance. By construction, $\phi_4^*$ conserves charge, as $B_4$ does:
the charges that either kernel assigns to the grid add up to the charge
of the particle wherever it lies, and both curves are zero. With the
$\psi_0$ kernel, the assigned charge varies with the position of the
particle, and its derivative jumps at $y_\alpha=0$, where the particle
crosses the boundary of a stencil.}
\label{fig:charge}
\end{figure}

The first of the desirable features of an SI kernel listed by Deserno
and Holm~\cite{DH1} is \emph{charge conservation}: wherever a particle
lies, the charges that spreading assigns to the grid points add up to
the charge of the particle. Spreading a charge $q$ at $\vy$ assigns to
the grid point $h\vell$ the charge $qh^3\Phi(h\vell-\vy)$, and since
$\Phi$ is a tensor product, these charges add up to
\begin{equation}
 \begin{gathered}
 qh^3\sum_{\vell\in\Z^3}\Phi(h\vell-\vy)=q\prod_{\alpha=1}^3S\Bigl(\frac{y_\alpha}h\Bigr),\\
 S\Bigl(\frac{y_\alpha}h\Bigr)=\sum_{\ell\in\Z}\phi\Bigl(\ell-\frac{y_\alpha}h\Bigr)
 =\sum_{m\in\Z}\widehat\phi(2\pi m)\,e^{2\pi\ii my_\alpha/h},
 \end{gathered}
 \label{eq:assigned-charge}
\end{equation}
where the second expression for $S$ follows from \cref{eq:poisson},
applied to $f(x)=\phi(x/h)$ with $d=1$ and $\xi=0$. The kernel thus
conserves charge if and only if $S(y_\alpha/h)=1$ for every $y_\alpha$,
that is, the integer translates $\phi(\cdot-\ell)$, $\ell\in\Z$, sum to
one. Since the Fourier
coefficients of $S$ are the values $\widehat\phi(2\pi m)$, and
$\widehat\phi(0)=1$, this holds if and only if $\widehat\phi$ vanishes at
every nonzero multiple of $2\pi$. The B-spline $B_P$ conserves charge,
since its transform vanishes to order $P$ at these points. The $\psi_0$
kernel does not: the total charge that it assigns varies with the
position of the particle relative to the grid (\cref{fig:charge}), and
\cref{sec:refinement} shows that, with AD, this variation produces a
force error that grows under refinement.

\subsection{The particle--mesh method}
\label{sec:pm-method}

Let $n_f=L/h$ be the number of grid points in each dimension, with $n_f$
even. The grid frequencies are $\vxi_\vk$ with $\vk$ in the index set
\begin{equation}
 \BZ(h)=\{-n_f/2,\ldots,n_f/2-1\}^3
 \label{eq:grid-indices}
\end{equation}
of the first Brillouin zone, so that each component of $\vxi_\vk$ is at
most the Nyquist frequency $\ximax(h)=\pi/h$ in magnitude. We retain a
set $\Modes(h)\subset\BZ(h)$ of nonzero indices with $|k_\alpha|<n_f/2$,
either all of them or those with $|\vxi_\vk|\le c_s/r_c$. We describe a
grid by the ratio of its Nyquist frequency to the band edge $c_s/r_c$,
\begin{equation}
 s=\frac{\ximax(h)}{c_s/r_c},\qquad\text{that is,}\qquad
 h=\frac{h_0}s,\quad h_0=\frac{\pi r_c}{c_s},
 \label{eq:grid-ratio}
\end{equation}
so that $h_0$ is the spacing whose Nyquist frequency is the band edge.
ESP's rule $n_f=\lceil c_sL/(\pi r_c)\rceil$~\cite{ESP} gives $s\ge1$,
close to one in large boxes; as the grid that ESP uses, we take the first
grid with even $n_f$ and $s\ge1$, and $s=1$ in the limit of large boxes. Upsampling, that is, refining the grid at fixed order
$P$, increases $s$. The
particle--mesh approximation of $u^s$ and
of the forces consists of five steps:
\begin{align}
 b_{\vell}&=\sum_{j=1}^Nq_j\,\Phi(h\vell-\vx_j),
 \label{eq:pm-spread}\\
 \widehat b_\vk&=h^3\sum_{\vell} b_{\vell}\,e^{\ii\vxi_\vk\cdot h\vell},
 \label{eq:pm-fft}\\
 \widehat c_\vk&=I(\vxi_\vk)\,\widehat b_\vk\quad(\vk\in\Modes(h)),
 \qquad \widehat c_\vk=0\quad\text{otherwise},
 \label{eq:pm-influence}\\
 c_{\vell}&=\frac1{L^3}\sum_{\vk\in\Modes(h)}\widehat c_\vk\,
          e^{-\ii\vxi_\vk\cdot h\vell},
 \label{eq:pm-ifft}\\
 u_{\rm mesh}(\vx)&=h^3\sum_{\vell} c_{\vell}\,\Phi(\vx-h\vell),\qquad
 \mathbf F_i^{\rm mesh}=-q_ih^3\sum_{\vell} c_{\vell}\,\nabla\Phi(\vx_i-h\vell).
 \label{eq:pm-gather}
\end{align}
Here $\vell$ runs over the $n_f^3$ grid points, $\Phi$ is extended
periodically, and $I$ is the influence function, a diagonal multiplier
in Fourier space. Steps \cref{eq:pm-fft} and \cref{eq:pm-ifft} are FFTs, which cost
$O(n_f^3\log n_f)$ operations; steps \cref{eq:pm-spread} and
\cref{eq:pm-gather} cost $O(NP^3)$.
The force in \cref{eq:pm-gather} is the AD force: it is the gradient of
the interpolated potential. Spectral differentiation, the ik method,
instead computes
the three components of the field,
$L^{-3}\sum_{\vk}\ii\xi_{\vk,\alpha}\widehat c_\vk
e^{-\ii\vxi_\vk\cdot h\vell}$, with three inverse FFTs and interpolates
each of them with $\Phi$.

To see what the method computes, fix $\vk\in\BZ(h)$ and denote the
aliases of $\vxi_\vk$ by
\begin{equation}
 \vxi_\vm=\vxi_\vk+\frac{2\pi}h\vm,\qquad \vm\in\Z^3;
 \label{eq:aliases}
\end{equation}
the alias with $\vm=\mathbf 0$ is $\vxi_\vk$ itself. We call the set of
these frequencies the class of $\vk$. Applying \cref{eq:poisson} to
\cref{eq:pm-spread}--\cref{eq:pm-fft} and using that $\Phi$ is even
gives
\begin{equation}
 \widehat b_\vk=\sum_{\vm\in\Z^3}\widehat\Phi(\vxi_\vm)\,\rhoh(\vxi_\vm),
 \label{eq:spread-aliases}
\end{equation}
so the grid coefficient of the frequency $\vxi_\vk$ collects $\rhoh$ at
$\vxi_\vk$ and at all its aliases, weighted by $\widehat\Phi$. In the same way,
\cref{eq:poisson} applied to \cref{eq:pm-gather} gives
\begin{equation}
 u_{\rm mesh}(\vx)=\frac1{L^3}\sum_{\vk\in\Modes(h)}I(\vxi_\vk)\,
   \widehat b_\vk
   \sum_{\vn\in\Z^3}\widehat\Phi(\vxi_\vn)\,e^{-\ii\vxi_\vn\cdot\vx},
 \label{eq:gather-aliases}
\end{equation}
so interpolation returns each grid coefficient at $\vxi_\vk$ and at all
its aliases. The main term, with $\vm=\vn=\mathbf 0$, is
$L^{-3}\sum_{\vk\in\Modes(h)}I(\vxi_\vk)\widehat\Phi(\vxi_\vk)^2
\rhoh(\vxi_\vk)e^{-\ii\vxi_\vk\cdot\vx}$, which equals the Fourier
series of $u^s$ truncated to $\Modes(h)$ when $I$ is the natural
influence function
\begin{equation}
 I_{\rm nat}(\vxi)=\frac{\Sh(\vxi)}{\widehat\Phi(\vxi)^2}.
 \label{eq:I-nat}
\end{equation}
ESP uses $I_{\rm nat}$. All other terms, with $\vm\ne\mathbf0$ or
$\vn\ne\mathbf0$, are the aliasing error. \Cref{sec:si} measures the aliasing error
and finds the influence function that minimizes it.

\section{Splitting kernels that are smooth at the cutoff}
\label{sec:split-design}

We turn now to the splitting kernel. Our goal is a split whose local
part is smooth at the cutoff, which pairs of particles cross many times
in a simulation, while the Fourier concentration of the smooth part stays
close to the optimum. In the splits of \cref{sec:splits}, the local part vanishes
beyond the cutoff, and its smoothness there is determined by the window
function alone. In this section, we state the smoothness that the time
integration requires and show that it amounts to linear conditions on the
window function. We then measure a split by its truncation error, show
that the $\psi_0$ split minimizes this error for the Coulomb interaction,
and construct the window function that satisfies the conditions with the
fewest and most concentrated PSWFs. The split depends only on the interaction, the cutoff, and the
bandwidth $c_s$. It fixes the Fourier modes that the particle--mesh
method must reproduce and their weights in the force error, and
\cref{sec:si} designs the SI kernel for them.

\subsection{Smoothness at the cutoff}
\label{sec:split-smoothness}

Molecular dynamics integrates Newton's equations of motion, most often by
the Verlet method, in its velocity or leap-frog form, which updates each
velocity by a quadrature rule of second order applied to the force along
the trajectory. Suppose that the energy is $d$
times continuously differentiable in the positions of the particles and
that its derivative of order $d+1$ jumps where a pair crosses the cutoff.
Then the jump adds a local error of order $d+1$ in the time step to the
step across that point, whereas the local error of the other steps is of
order three. For $d<2$, the jumps therefore dominate the errors of these
steps, and these errors can accumulate into a drift of the energy over many
time steps~\cite{CottrellTupper,EngleSkeel}. In simulations of a
Lennard-Jones fluid over $10^8$ time steps, for example, the energy
drifted steadily when the pair force jumped to zero at the cutoff, but it
was conserved when the force was shifted to vanish
there~\cite{ToxvaerdDyre} or the potential had four continuous
derivatives~\cite{ToxvaerdHeilmannDyre}. With $d=2$, every step keeps the
local order of the integrator. We therefore require the local part of
the split to be twice continuously differentiable at the cutoff; an
integrator of higher order would need more derivatives, and
\cref{sec:split-construction} gives splits with any number of them. In the
simulations of water in \cref{sec:nve}, the energy drifts steadily with
the kernels of ESP but stays close to its initial value when both
kernels are twice continuously differentiable.

In order to express this requirement in terms of the window function, we write
the local part of $t^{-\nu}$ in the construction of \cref{sec:splits} as
$t^{-\nu}-s_\nu(t)=t^{-\nu}R_\nu(t)$ for $0<t<1$, so that $R_\nu(t)$ is
the fraction of $t^{-\nu}$ that the local part carries; the local part
vanishes for $t>1$.
Since $t^{-\nu-2}=-(\nu t)^{-1}(t^{-\nu})'$, the local parts satisfy the
recurrence \cref{eq:inverse-recurrence} of the smooth parts, and
substituting $t^{-\nu}R_\nu$ into it and into \cref{eq:inverse-base} gives
\begin{equation}
 R_1(t)=\frac{\int_t^1w(\tau)\,d\tau}{\int_0^1w(\tau)\,d\tau},\qquad
 R_2=w,\qquad R_{\nu+2}=R_\nu-\frac t\nu R_\nu' .
 \label{eq:residual-recurrence}
\end{equation}
In particular,
\begin{equation}
 R_6(t)=w(t)-\tfrac58\,tw'(t)+\tfrac18\,t^2w''(t).
 \label{eq:R6}
\end{equation}
The local part of $r^{-\nu}$ is $\LL_\nu(r)=r^{-\nu}R_\nu(r/r_c)$ for
$0<r<r_c$, and for the Coulomb split, $R_1(t)=2\int_t^1\chi(\tau)\,d\tau$
by \cref{eq:coulomb-split}. The following proposition turns smoothness at
the cutoff into linear conditions on the window function.

\begin{proposition}
\label{prop:cutoff-matching}
Let $\nu\ge1$ and $d\ge0$, and let the window function $w$ be infinitely
differentiable on $[-1,1]$. The local part $\LL_\nu$ is $d$ times
continuously differentiable at $r=r_c$ if and only if
\begin{equation}
 R_\nu^{(j)}(1)=0,\qquad j=0,\ldots,d .
 \label{eq:cutoff-conditions}
\end{equation}
For $\nu=1$, these conditions hold if and only if $w^{(j)}(1)=0$ for
$j=0,\ldots,d-1$.
\end{proposition}

\begin{proof}
The local part and all its derivatives vanish for $r>r_c$, so $\LL_\nu$
is $d$ times continuously differentiable at $r_c$ if and only if its
derivatives of orders $0,\ldots,d$ tend to zero as $r\to r_c^-$, that is,
if and only if the derivatives of $t^{-\nu}R_\nu(t)$ of these orders
vanish at $t=1$. By the Leibniz rule, the $j$th derivative of
$t^{-\nu}R_\nu(t)$ at $t=1$ is $R_\nu^{(j)}(1)$ plus a linear combination
of $R_\nu(1),\ldots,R_\nu^{(j-1)}(1)$, so, by induction on $j$, these
conditions are equivalent to \cref{eq:cutoff-conditions}. For $\nu=1$, $R_1(1)=0$ and
$R_1^{(j)}(1)=-w^{(j-1)}(1)/\int_0^1w(\tau)\,d\tau$ for $j\ge1$.
\end{proof}

With $d=2$, the conditions are $w(1)=w'(1)=0$ for the Coulomb split, and
three linear conditions on the values of $w$ and its derivatives at $t=1$
for $\nu\ge2$. For $\nu=2$, they are $w(1)=w'(1)=w''(1)=0$, and they are
the same for $\nu=3$, since $R_3=R_1+tw/\int_0^1w(\tau)\,d\tau$ gives
$R_3(1)$, $R_3'(1)$, and $R_3''(1)$ as $w(1)$, $w'(1)$, and
$w'(1)+w''(1)$, respectively, divided by $\int_0^1w(\tau)\,d\tau$. For $\nu\ge4$, the
conditions also involve derivatives of order three and higher, and for
the dispersion interaction,
\cref{eq:R6} gives
\begin{equation}
 8w(1)-5w'(1)+w''(1)=0,\qquad 3w'(1)-3w''(1)+w'''(1)=0,\qquad
 w^{(4)}(1)=w'''(1).
 \label{eq:R6-conditions}
\end{equation}

The $\psi_0$ splits do not satisfy these conditions. For the Coulomb
split, the energy of a pair is continuous at the cutoff, but
$R_1'(1)=-2\chi(1)$, so the force between the two particles jumps there by
$2\chi(1)$ times the bare force $q_1q_2/r_c^2$; by ESP's rule
$\chi(1)=\varepsilon$, this factor is twice the tolerance. For $\nu\ge2$, the local
part itself jumps, so the energy of a pair changes by $R_\nu(1)$ times the
bare energy $q_1q_2r_c^{-\nu}$ when the pair crosses the cutoff. For
$\nu=2$ and $\nu=3$, $R_\nu(1)$ is $w(1)$ and $2\chi(1)$, respectively,
of the order of the tolerance. For $\nu\ge4$, $R_\nu(1)$ contains derivatives of $w$ at
$t=1$, which the prolate equation makes large: by \cref{eq:prolate-ode},
$w'(1)=\frac12(\lambda_0-c_s^2)w(1)$ with $\lambda_0\le c_s^2/3$, and
$w''(1)$ is larger still. For the dispersion interaction, $R_6(1)$ is
$36$ to $1810$ times the tolerance for the tolerances from $10^{-2}$ to
$10^{-8}$.

\subsection{The truncation error}
\label{sec:split-truncation}

The grid that ESP uses retains a cube of modes that contains the band
$|\vxi|\le c_s/r_c$, and we measure a split by the force that the modes
outside the band carry.
Let $\Sh_\nu$ be the Fourier transform of the smooth part $\Sp_\nu$ of
$r^{-\nu}$, so that $\Sh_1=\Sh$. Omitting the modes outside the band from
the smooth potential \cref{eq:split-potential} of a unit charge at $\vy$
changes the force on a unit charge at $\vx$ by the gradient of the
omitted terms. When $\vx$ and $\vy$ are independent and uniformly
distributed in the box, Parseval's identity gives the mean square of this
change as $L^{-6}\sum|\vxi_\vk|^2\Sh_\nu(\vxi_\vk)^2$, summed over the
omitted modes. As $L\to\infty$, $L^{-3}$ times such sums tend to integrals
over $\R^3$ divided by $(2\pi)^3$, and we measure a split by its relative
truncation error
\begin{equation}
 \biggl(\frac{\int_{|\vxi|>c_s/r_c}|\vxi|^2\Sh_\nu(\vxi)^2\,d\vxi}
       {\int_{\R^3}|\vxi|^2\Sh_\nu(\vxi)^2\,d\vxi}\biggr)^{1/2},
 \label{eq:truncation-error}
\end{equation}
the RMS force of the omitted modes relative to the RMS force of the
smooth part. It is defined only when $R_\nu(1)=0$: for the $\psi_0$ splits
with $\nu\ge2$, the smooth part jumps at the cutoff, and both mean-square
forces are infinite. For the Coulomb interaction, the relative truncation
error depends only on the Fourier transform of the splitting kernel, and
the $\psi_0$ split minimizes it.

\begin{lemma}
\label{lem:psi0-optimal}
Let $\chi$ be a real, even, square-integrable function on $[-1,1]$ with
$\int_{-1}^1\chi(t)\,dt=1$, and let $\Sh$ be the Fourier transform
\cref{eq:split-transform} of the smooth part of the Coulomb split
\cref{eq:coulomb-split} with the cutoff $r_c$. Then, for every $c_s>0$,
\begin{equation}
 \frac{\int_{|\vxi|>c_s/r_c}|\vxi|^2\Sh(\vxi)^2\,d\vxi}
      {\int_{\R^3}|\vxi|^2\Sh(\vxi)^2\,d\vxi}
 =\frac{\int_{|\omega|>c_s}\widehat\chi(\omega)^2\,d\omega}
       {\int_{\R}\widehat\chi(\omega)^2\,d\omega}.
 \label{eq:tail-band-energy}
\end{equation}
Among all such $\chi$, the $\psi_0$ split \cref{eq:psi0-split} alone
minimizes the relative truncation error \cref{eq:truncation-error}.
\end{lemma}

\begin{proof}
By \cref{eq:split-transform},
$|\vxi|^2\Sh(\vxi)^2=16\pi^2\widehat\chi(|\vxi|r_c)^2/|\vxi|^2$. In
spherical coordinates, the radial measure $4\pi|\vxi|^2\,d|\vxi|$ cancels
the factor $|\vxi|^{-2}$, and the substitution $\omega=|\vxi|r_c$ turns
the numerator and the denominator on the left of
\cref{eq:tail-band-energy} into $64\pi^3/r_c$ times the integrals of
$\widehat\chi(\omega)^2$ over $\omega>c_s$ and over $\omega>0$,
respectively; since
$\widehat\chi$ is even, this gives \cref{eq:tail-band-energy}. Its right
side is one minus the fraction of the energy of $\widehat\chi$ in
$[-c_s,c_s]$, and among the functions supported on $[-1,1]$, this fraction
is largest exactly for the multiples of $\psi_0^{c_s}$ (\cref{sec:pswf}).
With $\int_{-1}^1\chi(t)\,dt=1$, the minimizer is
$\chi=\psi_0^{c_s}/\int_{-1}^1\psi_0^{c_s}(t)\,dt$.
\end{proof}

In other words, the relative truncation error of a Coulomb split is the
square root of the fraction of the energy of $\widehat\chi$ outside
$[-c_s,c_s]$, so concentration in frequency, the property for which the
DMK method~\cite{DMK} and ESP~\cite{ESP} chose $\psi_0$, is exactly the
right measure for the Coulomb split. For $\chi=\sum_na_n\psi_n^{c_s}$,
\cref{eq:pswf-out-of-band} makes this fraction the average of the
corresponding fractions of the $\psi_n^{c_s}$, weighted by $a_n^2$, and
these grow with $n$.

\subsection{The \texorpdfstring{$C^2$}{C2} splits}
\label{sec:split-construction}

For the Coulomb split, the window function that minimizes the truncation error
subject to the conditions of \cref{prop:cutoff-matching} does not exist.
The fraction \cref{eq:tail-band-energy} is continuous in $L^2(-1,1)$, but the values of a window function and its derivatives at $t=1$
are not. Multiplying $\psi_0^{c_s}$ by a smooth function that equals one
outside layers of width $\delta$ at $\pm1$ and vanishes to second order at
$\pm1$ gives window functions that satisfy the conditions and tend to
$\psi_0^{c_s}$ in $L^2$ as $\delta\to0$. Their truncation errors approach
that of the $\psi_0$ split, which no window function that satisfies the
conditions attains. But they confine the change of the split to a layer
at the cutoff that a pair crosses within a single time step, so that
the integrator still sees a jump.

We therefore take the window function in the span of the fewest even PSWFs that
admit the conditions, so that its change from the $\psi_0$ split extends
over the whole local part rather than a thin layer. Requiring $d$
continuous derivatives of the local part at the cutoff imposes $m$
conditions, $m=d$
for $\nu=1$ and $m=d+1$ for $\nu\ge2$, and the split has the window
function
\begin{equation}
 w(t)=\sum_{i=0}^{m}w_{2i}\,\psi_{2i}^{c_s}(t),\qquad |t|\le1,
 \label{eq:split-expansion}
\end{equation}
with a bandwidth $c_s$ and the normalization $w(0)=1$. Without conditions, this is the $\psi_0$ split, and since the fractions
of the $\psi_n^{c_s}$ outside the band grow with $n$, the split uses
the most concentrated PSWFs that admit the conditions. The $C^2$ split,
$d=2$, has $m=2$ for $\nu=1$ and $m=3$ for $\nu\ge2$. When the conditions
are $w^{(j)}(1)=0$ for $j<m$, as for $\nu\le3$, the coefficients have a
closed form.

\begin{theorem}
\label{thm:split-closed-form}
Let $c>0$ and $m\ge1$. There is exactly one function
$w=\sum_{i=0}^mw_{2i}\psi_{2i}^c$ with $w(0)=1$ and $w^{(j)}(1)=0$ for
$j=0,\ldots,m-1$. Its coefficients are positive and proportional to
\begin{equation}
 \frac1{\psi_{2i}^c(1)\prod_{k\ne i}(\lambda_{2i}-\lambda_{2k})},
 \qquad i=0,\ldots,m,
 \label{eq:split-closed-form}
\end{equation}
where the product runs over $k=0,\ldots,m$.
\end{theorem}

\begin{proof}
The prolate equation \cref{eq:prolate-ode} holds at $t=1$ by continuity,
since $\psi_n$ is entire, and there the coefficient of $\psi''$ vanishes,
so the equation gives
$2\psi_n'(1)=(\lambda_n-c^2)\psi_n(1)$, where $\psi_n=\psi_n^c$.
Differentiating the prolate equation $j$ times and setting $t=1$ gives
\[
 2(j+1)\psi_n^{(j+1)}(1)=\bigl(\lambda_n-c^2-j(j+1)\bigr)\psi_n^{(j)}(1)
   -2jc^2\psi_n^{(j-1)}(1)-j(j-1)c^2\psi_n^{(j-2)}(1).
\]
By induction on $j$, $\psi_n^{(j)}(1)=p_j(\lambda_n)\psi_n(1)$, where
$p_j$ is a polynomial of degree $j$ whose coefficients depend only on $c$
and whose leading coefficient is $1/(2^jj!)$. In particular,
$\psi_n(1)\ne0$, since otherwise all derivatives of the analytic
function $\psi_n$ would vanish at $t=1$. With $u_i=w_{2i}\psi_{2i}(1)$, the
conditions read $\sum_iu_ip_j(\lambda_{2i})=0$ for $j<m$, and since
$p_0,\ldots,p_{m-1}$ span the polynomials of degree less than $m$, the
conditions are equivalent to $\sum_iu_i\lambda_{2i}^j=0$ for $j<m$. The eigenvalues
are distinct, so this Vandermonde system of $m$ equations in $m+1$
unknowns has a one-dimensional space of solutions, spanned by the weights
$u_i=1/\prod_{k\ne i}(\lambda_{2i}-\lambda_{2k})$ of the divided difference
of order $m$, which vanishes on the polynomials of degree less than $m$.
This gives \cref{eq:split-closed-form}. Since
$\lambda_0<\lambda_2<\cdots<\lambda_{2m}$, the product in
\cref{eq:split-closed-form} has the sign $(-1)^{m-i}$, and since
$\psi_{2i}$ is positive at $t=0$ and has exactly $i$ simple zeros in
$(0,1)$, $\psi_{2i}(1)$ has the sign $(-1)^i$. The coefficients of every
solution therefore have one sign, and since every $\psi_{2i}(0)$ is
positive, $w(0)$ has that sign as well; the normalization $w(0)=1$ makes
every coefficient positive.
\end{proof}

In other words, the window function is a divided difference of order $m$ of the
functions $\psi_{2i}^c/\psi_{2i}^c(1)$ with respect to their eigenvalues
$\lambda_{2i}$, scaled so that $w(0)=1$. The $C^2$ Coulomb split uses
$\psi_0^{c_s}$, $\psi_2^{c_s}$, and $\psi_4^{c_s}$, and the $C^2$ splits of
$1/r^2$ and $1/r^3$ share one window function in the span of
$\psi_0^{c_s},\ldots,\psi_6^{c_s}$. For $\nu\ge4$, the recursion in the
proof turns the three conditions into three linear equations for the four
coefficients, in which the entries multiplying $w_{2i}$ are polynomials
in $\lambda_{2i}$ times $\psi_{2i}(1)$. We compute the coefficients as the
null vector of this $3\times4$ system by the singular value
decomposition. At every bandwidth
used in \cref{sec:numerics}, the system has rank three and the solution
has $w(0)\ne0$, so the split is unique. For every $\nu$, the
construction thus requires no minimization, and \cref{sec:split-results}
gives the truncation errors of the $C^2$ splits.

\section{SI kernels that minimize the force error}
\label{sec:si}

We turn now to the SI kernel. With the split of \cref{sec:split-design}
fixed, the particle--mesh method of \cref{sec:pm-method} approximates the
smooth part, and the SI kernel determines the aliasing error of this
approximation. In order to design the kernel, we first express the RMS
force error in terms of the Fourier transforms of the two kernels
(\cref{sec:pair-self,sec:pair-fourier,sec:alias-tail}). The resulting
formula shows that with AD, the error of an SI kernel of fixed shape
tends to zero under refinement of the grid only if the kernel conserves
charge and, if it is continuous, also the dipole, and that without
charge conservation, the error grows with the size of the box at a fixed
grid spacing (\cref{sec:refinement}). We then construct, for each order and
tolerance, the kernel with these properties that meets the tolerance on
the coarsest grid (\cref{sec:si-design,sec:si-construction}), and we end
with the choice of the parameters (\cref{sec:parameters}).

\subsection{Pair and self terms}
\label{sec:pair-self}

Since the local part is summed exactly, the error of the computed force
on particle $i$ is that of its smooth part,
$\mathbf F_i^{\rm mesh}+q_i\nabla u^s(\vx_i)$. We measure accuracy by the
RMS force error, the square root of the mean-square force error
\begin{equation}
 \mathcal V_F=\frac1N\,\mathbb E\sum_{i=1}^N
              \bigl|\mathbf F_i^{\rm mesh}+q_i\nabla u^s(\vx_i)\bigr|^2 ,
 \label{eq:rms-force}
\end{equation}
where $\mathbb E$ denotes the mean over configurations; for independent,
uniformly distributed positions, the RMS force error is a normalized
$L^2$ norm of the force error over the positions of the particles. Let $e(\vx,\vy)$ be the
computed minus the exact smooth potential at $\vx$ produced by a unit
charge at $\vy$ and its periodic images, and $e(\vx,\vx)$ the error of
the potential that the charge produces at its own position, after the
subtraction of $\Sp(0)$ in \cref{eq:split-potential}. Spreading and
interpolation use the same kernel and $I$ is real and even, so
$e(\vx,\vy)=e(\vy,\vx)$, and with AD
\begin{equation}
 \mathbf F_i^{\rm mesh}+q_i\nabla u^s(\vx_i)
 =-\sum_{j\ne i}q_iq_j\,\nabla_\vx e(\vx_i,\vx_j)
  -\frac12q_i^2\,\nabla[e(\vx,\vx)]_{\vx=\vx_i},
 \label{eq:force-error-pairs}
\end{equation}
where $\nabla_\vx$ acts on the first argument. The last term is the self
force, which depends on the position of the particle relative to the
grid~\cite{SelfP3M,Ball}. We define the mean-square pair and self force
errors
\begin{equation}
 \begin{gathered}
 \sigma_{\rm pair,F}^2=\frac1{L^6}\int_{[0,L)^3}\int_{[0,L)^3}
      |\nabla_\vx e(\vx,\vy)|^2\,d\vx\,d\vy,\\
 \sigma_{\rm self,F}^2=\frac1{L^3}\int_{[0,L)^3}
      \Bigl|\tfrac12\nabla[e(\vx,\vx)]\Bigr|^2\,d\vx .
 \end{gathered}
 \label{eq:pair-self-errors}
\end{equation}
Thus $\sigma_{\rm pair,F}^2$ is the mean-square error of the force
between two unit charges at independent, uniformly distributed
positions, which is the error measure of Hockney and Eastwood~\cite{HE}
divided by $L^3$; $\sigma_{\rm self,F}$ is the RMS self force of a unit
charge. For the $\psi_0$ kernel, which AD codes differentiate only inside
its support (\cref{sec:alias-tail}), we use the forces that these codes
compute. For random charges, the standard estimate of the RMS force
error~\cite{Ball,DH2} is computed from $\sigma_{\rm pair,F}$; the following
proposition gives the mean-square force error exactly, including the
self term.

\begin{proposition}
\label{prop:force-variance}
Suppose that the positions of the particles are independent and
uniformly distributed in $[0,L)^3$, that the charges are fixed, and that
\begin{equation}
 \int_{[0,L)^3}\nabla_\vx e(\vx,\vy)\,d\vy=\mathbf 0
 \qquad\text{for every }\vx .
 \label{eq:zero-mean-pair}
\end{equation}
Then
\begin{equation}
 \mathcal V_F=\frac{(M_2^2-M_4)\,\sigma_{\rm pair,F}^2
                 +M_4\,\sigma_{\rm self,F}^2}{N},
 \qquad M_2=\sum_{i=1}^Nq_i^2,\quad M_4=\sum_{i=1}^Nq_i^4 .
 \label{eq:force-variance}
\end{equation}
The particle--mesh method of \cref{sec:pm-method} satisfies
\cref{eq:zero-mean-pair}.
\end{proposition}

\begin{proof}
Expand the square of \cref{eq:force-error-pairs} for the random positions.
For $j\ne k$, both different from $i$, $\vx_j$ and $\vx_k$ are independent
for fixed $\vx_i$, so, by \cref{eq:zero-mean-pair}, the product of the terms
for $j$ and $k$ has mean zero, and so does the product of a pair term with
the self term, which depends only on $\vx_i$. The mean square of the
$i$th force error is therefore
$\sum_{j\ne i}q_i^2q_j^2\sigma_{\rm pair,F}^2+q_i^4\sigma_{\rm self,F}^2$,
and summing over $i$ with $\sum_i\sum_{j\ne i}q_i^2q_j^2=M_2^2-M_4$ gives
\cref{eq:force-variance}. For \cref{eq:zero-mean-pair}, averaging over
the position of the unit source in the box replaces the source by the
uniform density $L^{-3}$: by \cref{eq:pm-spread}, the grid then holds the constant
$L^{-3}$, whose only Fourier mode, the zero mode, is not in $\Modes(h)$;
the exact smooth potential contains only nonzero modes
$e^{\ii\vxi_\vk\cdot\vy}$, each of which averages to zero over $\vy$.
\end{proof}

For unit charges, $M_2=M_4=N$ and
$\mathcal V_F=(N-1)\sigma_{\rm pair,F}^2+\sigma_{\rm self,F}^2$. Since the
pair term grows with the number of charges in a given box and the self
term does not, P3M is designed for the pair term alone~\cite{EnergyP3M,HE},
and so is the SI kernel below; this keeps the kernel independent of the
charges of the system. The self force of each particle can be computed
from the mesh potential of the particle's own charge and subtracted from
the force on the particle~\cite{SelfP3M}, and \cref{sec:complete-results}
compares the two terms in complete calculations.

\subsection{The pair error in Fourier space}
\label{sec:pair-fourier}

Hockney and Eastwood showed that the pair error is a quadratic in the
values of the influence function at the grid frequencies~\cite{HE}. With
the aliases \cref{eq:aliases} and AD, it reads~\cite{Ball,Stern}
\begin{equation}
 \sigma_{\rm pair,F}^2=\frac1{L^6}\sum_{\vk\in\BZ(h)}
  \bigl(A_\vk I(\vxi_\vk)^2-2B_\vk I(\vxi_\vk)+C_\vk\bigr),
 \label{eq:pair-quadratic}
\end{equation}
with $I(\vxi_\vk)=0$ for $\vk\notin\Modes(h)$, where
\begin{equation}
 \begin{gathered}
 A_\vk=\sum_\vm\widehat\Phi(\vxi_\vm)^2
       \sum_\vn|\vxi_\vn|^2\widehat\Phi(\vxi_\vn)^2,\\
 B_\vk=\sum_\vm|\vxi_\vm|^2\widehat\Phi(\vxi_\vm)^2\Sh(\vxi_\vm),\qquad
 C_\vk=\sum_\vm|\vxi_\vm|^2\Sh(\vxi_\vm)^2,
 \end{gathered}
 \label{eq:ABC}
\end{equation}
and the term with $\vm=\mathbf 0$ is omitted for $\vk=\mathbf 0$. Here the
source aliases $\vxi_\vm$ of the class of $\vk$ enter through spreading
\cref{eq:spread-aliases} and the target aliases $\vxi_\vn$ through
interpolation \cref{eq:gather-aliases}. Spectral differentiation
multiplies the grid coefficient by $\ii\vxi_\vk$ before interpolation, and
then $|\vxi_\vn|^2$ in $A_\vk$ becomes $|\vxi_\vk|^2$ and $|\vxi_\vm|^2$ in
$B_\vk$ becomes $\vxi_\vk\cdot\vxi_\vm$. For $A_\vk>0$, the pair error is
minimized by the optimal influence function
\begin{equation}
 I_{\rm opt}(\vxi_\vk)=\frac{B_\vk}{A_\vk},
 \qquad\text{with minimum}\quad
 \frac1{L^6}\sum_{\vk\in\Modes(h)}\Bigl(C_\vk-\frac{B_\vk^2}{A_\vk}\Bigr)
 +\frac1{L^6}\sum_{\vk\in\BZ(h)\setminus\Modes(h)}C_\vk ,
 \label{eq:I-opt}
\end{equation}
which is the influence function of Hockney and Eastwood~\cite{HE} for
spectral differentiation and that of~\cite{Ball,Stern} for AD. With
$I_{\rm opt}$, the pair error depends only on the two kernels and the
grid, and we use $I_{\rm opt}$ in all comparisons of SI kernels for the
Coulomb interaction.

\subsection{The aliasing error and the truncation error}
\label{sec:alias-tail}

The SI kernel enters \cref{eq:pair-quadratic} only through the retained
classes $\vk\in\Modes(h)$. To isolate its effect, we compare the mesh
result with the Fourier series of $u^s$ truncated to $\Modes(h)$, which
keeps in each retained class only the frequency $\vxi_\vk$ itself. The
difference between the full and the truncated reference is the
truncation error of the split. Its mean-square force,
$L^{-6}\sum|\vxi|^2\Sh(\vxi)^2$ over the omitted frequencies, depends
only on the split and on $\Modes(h)$. When $\Modes(h)$ consists of the
modes with $|\vxi_\vk|\le c_s/r_c$, the ratio of this mean-square force to
that of the smooth part tends, as $L\to\infty$, to the square of the relative
truncation error \cref{eq:truncation-error}. The complete pair error also
contains a cross term between the aliasing error and the truncation
error, since the target aliases of a retained class lie among the omitted
frequencies; the complete calculations of \cref{sec:complete-results}
include all three parts.

With the truncated reference, the sums for $B_\vk$ and $C_\vk$ in
\cref{eq:ABC} reduce to their terms with $\vm=\mathbf 0$:
$B_\vk=\widehat\Phi(\vxi_\vk)^2|\vxi_\vk|^2\Sh(\vxi_\vk)$ and
$C_\vk=|\vxi_\vk|^2\Sh(\vxi_\vk)^2$ for $\vk\in\Modes(h)$. Put
\begin{equation}
 T(\vk)=\sum_{\vm}\frac{\widehat\Phi(\vxi_\vm)^2}{\widehat\Phi(\vxi_\vk)^2},
 \qquad
 T_2(\vk)=\sum_{\vm}|\vxi_\vm|^2\,
         \frac{\widehat\Phi(\vxi_\vm)^2}{\widehat\Phi(\vxi_\vk)^2}.
 \label{eq:alias-ratios}
\end{equation}
Then $A_\vk=\widehat\Phi(\vxi_\vk)^4T T_2$ for AD and
$A_\vk=\widehat\Phi(\vxi_\vk)^4|\vxi_\vk|^2T^2$ for spectral
differentiation, and \cref{eq:I-opt}, with $B_\vk$ for the truncated
reference, gives the aliasing errors
\begin{equation}
 \begin{gathered}
 Q_{\rm AD}=\sum_{\vk\in\Modes(h)}|\vxi_\vk|^2\Sh(\vxi_\vk)^2\,\eta_{\rm AD}(\vk),\qquad
 Q_{\rm ik}=\sum_{\vk\in\Modes(h)}|\vxi_\vk|^2\Sh(\vxi_\vk)^2\,\eta_{\rm ik}(\vk),\\
 \eta_{\rm AD}(\vk)=1-\frac{|\vxi_\vk|^2}{T(\vk)T_2(\vk)},\qquad
 \eta_{\rm ik}(\vk)=1-T(\vk)^{-2},
 \end{gathered}
 \label{eq:Q-AD-ik}
\end{equation}
so that $\sigma_{\rm pair,F}^2=L^{-6}Q_{\rm AD}$ for AD against the
truncated reference with $I=I_{\rm opt}$. In the truncated reference,
the mode $\vxi_\vk$ contributes $L^{-6}|\vxi_\vk|^2\Sh(\vxi_\vk)^2$ to
the mean-square force, so $\eta_{\rm AD}(\vk)$ and $\eta_{\rm ik}(\vk)$
are the ratios of the mean-square aliasing error of this mode to its
mean-square force. Since $|\vxi_\vm|\ge|\vxi_\vk|$ for every $\vm$,
$T_2\ge|\vxi_\vk|^2T$, and hence $\eta_{\rm AD}(\vk)\ge\eta_{\rm ik}(\vk)$:
with the same continuous SI kernel, AD is never more accurate than
spectral differentiation. The difference is the weight $|\vxi_\vm|^2/|\vxi_\vk|^2$
of each alias ratio $\widehat\Phi(\vxi_\vm)^2/\widehat\Phi(\vxi_\vk)^2$ in
$T_2$. This weight is large for the modes with $|\vxi_\vk|\ll2\pi/h$ that
carry most of the Coulomb force: AD differentiates the target aliases that
interpolation adds, which oscillate on the scale of the grid.

For a tensor-product kernel \cref{eq:si-kernel},
\begin{equation}
 T(\vk)=\prod_{\alpha=1}^3\bigl(1+\mathcal A(h\xi_{\vk,\alpha})\bigr),
 \qquad
 \mathcal A(q)=\sum_{m\ne0}\left|\frac{\widehat\phi(q+2\pi m)}{\widehat\phi(q)}\right|^2,
 \label{eq:alias-power}
\end{equation}
and
\begin{equation}
 T_2(\vk)=\sum_{\alpha=1}^3\Bigl[\sum_{m\in\Z}
   \Bigl(\xi_{\vk,\alpha}+\frac{2\pi m}h\Bigr)^2
   \left|\frac{\widehat\phi(h\xi_{\vk,\alpha}+2\pi m)}{\widehat\phi(h\xi_{\vk,\alpha})}\right|^2\Bigr]
   \prod_{\beta\ne\alpha}\bigl(1+\mathcal A(h\xi_{\vk,\beta})\bigr).
 \label{eq:T2-tensor}
\end{equation}
The one-dimensional alias power $\mathcal A(q)$ is an infinite
sum, but it equals a mean over a grid cell of a finite sum. By
\cref{eq:poisson} and Parseval's identity in $y_\alpha$,
\begin{equation}
 \begin{aligned}
 |\widehat\phi(q)|^2\bigl(1+\mathcal A(q)\bigr)&=\sum_{m\in\Z}|\widehat\phi(q+2\pi m)|^2\\
 &=\frac1h\int_0^h\Bigl|\sum_{j\in\Z}\phi(j-y_\alpha/h)\,e^{\ii q(j-y_\alpha/h)}\Bigr|^2dy_\alpha ,
 \end{aligned}
 \label{eq:alias-cell}
\end{equation}
where the inner sum has at most $P+1$ nonzero terms. The sums in
\cref{eq:T2-tensor} follow in the same way, with $\phi$ replaced by $\phi'$. We evaluate these
cell integrals by Gauss quadrature on the intervals on which the
integrand is smooth, and so account for all aliases.

The $\psi_0$ kernel requires one more step. Since it jumps at the ends of
its support, integration by parts shows that $\widehat\phi(q)$ decays only like
$1/|q|$, so the terms $|\vxi_\vm|^2\widehat\Phi(\vxi_\vm)^2$ do not tend
to zero and $T_2$ is infinite. AD codes therefore differentiate this
kernel only inside its support~\cite{ESP}. We call this derivative,
extended by zero outside the support, the interior derivative and denote
it by $\phi'$; for a continuous kernel, it is the derivative on $\R$. If $\phi$ is even and
continuous inside its support, as every kernel in this paper is,
integration by parts gives
\begin{equation}
 \widehat{\phi'}(q)=-\ii q\widehat\phi(q)+2\ii\phi(P/2^-)\sin(qP/2).
 \label{eq:interior-derivative}
\end{equation}
The last term vanishes for a continuous kernel. For the $\psi_0$ kernel,
it distinguishes $\phi'$ from the derivative of $\phi$ in the sense of
distributions, whose Fourier transform is $-\ii q\widehat\phi(q)$. We
treat the $\psi_0$ kernel as AD codes do: in the derivation of
\cref{eq:pair-quadratic}, the Fourier transform of $\nabla\Phi$ at the
target alias $\vxi_\vn$ is replaced by that of the gradient taken inside
the support. In
\cref{eq:T2-tensor}, the replacement turns
$(\xi_{\vk,\alpha}+2\pi m/h)^2|\widehat\phi(h\xi_{\vk,\alpha}+2\pi m)|^2$ into
$h^{-2}|\widehat{\phi'}(h\xi_{\vk,\alpha}+2\pi m)|^2$; by
\cref{eq:interior-derivative}, the two agree for a continuous kernel. The
replacement also changes the main term $\vm=\mathbf 0$ of $B_\vk$, and with
the truncated reference, $\eta_{\rm AD}$ in \cref{eq:Q-AD-ik} becomes
\begin{equation}
 \eta_{\rm AD}(\vk)=1-\frac{|\vxi_\vk|^2}{T(\vk)T_2(\vk)}
  \Bigl(1-\frac{2\phi(P/2^-)}{h|\vxi_\vk|^2}\sum_{\alpha=1}^3
  \frac{\xi_{\vk,\alpha}\sin(Ph\xi_{\vk,\alpha}/2)}{\widehat\phi(h\xi_{\vk,\alpha})}\Bigr)^2,
 \label{eq:eta-AD-jump}
\end{equation}
which reduces to \cref{eq:Q-AD-ik} for a continuous kernel. The
interior derivative is not the gradient of the computed energy, which
jumps when a particle crosses the boundary of a stencil.

\subsection{Upsampling and the size of the box}
\label{sec:refinement}

The error formulas just derived explain why, with AD, upsampling reduces
the error of the $\psi_0$ kernel far less than that of B-splines, and why,
at a fixed grid spacing, its error grows with the size of the box. They also show which SI kernels
avoid both defects. Consider upsampling, $s\to\infty$ in
\cref{eq:grid-ratio}, with the order and the kernel $\phi$ held fixed. A
component $\xi$ of a fixed grid frequency
$\vxi_\vk$ then enters \cref{eq:alias-power,eq:T2-tensor} through
$q=h\xi$, which tends to zero as $s\to\infty$. This limit is not a
practical regime, but it explains the trend at the moderate upsampling
factors used in practice (\cref{sec:band-results}). The same variable is small for the lowest
modes of a large box on a fixed grid, since a box of side $L$ has modes
with $|\vxi_\vk|=2\pi/L$. Recall that $S$ is the one-dimensional
factor \cref{eq:assigned-charge} of the assigned charge. The following
proposition covers B-splines, whose derivatives jump at their knots, and
the $\psi_0$ kernel, which jumps at the ends of its support.

\begin{proposition}
\label{prop:alias-limit}
Let $\phi$ be even, supported in $[-P/2,P/2]$ for an integer $P\ge2$,
and normalized by $\int_{\R}\phi(t)\,dt=1$. Let its restriction to
$[-P/2,P/2]$ be continuous and piecewise $C^2$, so that $\phi$ can jump
only at $\pm P/2$. Let $\phi'$ be its derivative inside the support,
extended by zero. Then, for every $\xi\in\R$,
\begin{align}
 \lim_{h\to0}\mathcal A(h\xi)&=\sum_{m\ne0}|\widehat\phi(2\pi m)|^2,
 \label{eq:alias-limit}\\
 \lim_{h\to0}\sum_{m\ne0}\frac{|\widehat{\phi'}(h\xi+2\pi m)|^2}{|\widehat\phi(h\xi)|^2}
   &=\sum_{m\ne0}(2\pi m)^2|\widehat\phi(2\pi m)|^2 .
 \label{eq:alias-limit-ad}
\end{align}
Both limits vanish if and only if the integer translates of $\phi$ sum
to one. If they do not, then for every $\vk\in\Z^3\setminus\{\mathbf 0\}$,
\begin{equation}
 \lim_{h\to0}\eta_{\rm ik}(\vk)
   =1-\Bigl(1+\sum_{m\ne0}|\widehat\phi(2\pi m)|^2\Bigr)^{-6}>0,
 \qquad
 \lim_{h\to0}\eta_{\rm AD}(\vk)=1 .
 \label{eq:ratio-limit}
\end{equation}
\end{proposition}

\begin{proof}
By \cref{eq:interior-derivative},
$\widehat{\phi'}(2\pi m)=-2\pi\ii m\widehat\phi(2\pi m)$, since $P$ is an
integer. Integration by parts on each interval on which $\phi$ is $C^2$
gives $|\widehat\phi(q)|+|\widehat{\phi'}(q)|\le C/(1+|q|)$ for some
constant $C$. For $|q|\le\pi$, the terms of the sums over $m\ne0$ in
\cref{eq:alias-limit,eq:alias-limit-ad} are therefore bounded by a
summable sequence in $m$, and since $\widehat\phi(0)=1$, dominated
convergence gives \cref{eq:alias-limit,eq:alias-limit-ad}. Either limit
vanishes exactly when $\widehat\phi(2\pi m)=0$ for every $m\ne0$, that
is, by \cref{sec:si-kernels}, when the integer translates of $\phi$ sum
to one.

Next, \cref{eq:alias-power,eq:alias-limit} give
$\lim_{h\to0}T(\vk)=(1+\sum_{m\ne0}|\widehat\phi(2\pi m)|^2)^3$, and
hence the limit of $\eta_{\rm ik}(\vk)$ in \cref{eq:ratio-limit}. With AD, by
\cref{eq:interior-derivative} and the replacement described in
\cref{sec:alias-tail}, the bracket in \cref{eq:T2-tensor} equals
$h^{-2}\sum_{m\in\Z}|\widehat{\phi'}(h\xi_{\vk,\alpha}+2\pi m)|^2/|\widehat\phi(h\xi_{\vk,\alpha})|^2$.
After multiplication by $h^2$, its term with $m=0$ tends to
$|\widehat{\phi'}(0)|^2=0$, so
\cref{eq:alias-limit,eq:alias-limit-ad} give
\[
 \lim_{h\to0}h^2T_2(\vk)
 =3\Bigl(1+\sum_{m\ne0}|\widehat\phi(2\pi m)|^2\Bigr)^2\sum_{m\ne0}(2\pi m)^2|\widehat\phi(2\pi m)|^2,
\]
which is positive, and hence $\lim_{h\to0}|\vxi_\vk|^2/(T(\vk)T_2(\vk))=0$.
For a kernel that jumps, the squared factor in \cref{eq:eta-AD-jump}
tends to $(1-P\phi(P/2^-))^2$, so $\eta_{\rm AD}(\vk)\to1$ in either case.
\end{proof}

With spectral differentiation, a kernel whose integer translates do not
sum to one thus leaves an error floor at each mode. With AD, the
aliasing error of each mode becomes as large as the force of the mode,
for the following reason. Interpolating a constant mesh potential $u_0$ with $\nabla\Phi$ gives a
unit charge at $\vy$ the force
$-u_0\nabla_\vy\bigl(h^3\sum_{\vell}\Phi(h\vell-\vy)\bigr)$, proportional
to the gradient of the charge that it assigns to the grid. By \cref{eq:assigned-charge},
the $\alpha$th component of this force is
$-u_0h^{-1}S'(y_\alpha/h)\prod_{\beta\ne\alpha}S(y_\beta/h)$, so the force
grows like $1/h$; by Parseval's identity, the mean square of
$S'(y_\alpha/h)$ over a grid cell is the limit in \cref{eq:alias-limit-ad}.

The $\psi_0$ kernel has this defect whatever its bandwidth.

\begin{proposition}
\label{prop:psi0-defect}
For every $c_w>0$, the $\psi_0$ kernel \cref{eq:psi0-kernel} satisfies
\begin{equation}
 \sum_{m\ne0}|\widehat\phi(2\pi m)|^2>0 .
 \label{eq:psi0-defect}
\end{equation}
On every interval
$0<c_-\le c_w\le c_+<\infty$, this quantity, as a function of $c_w$, has
a positive minimum.
\end{proposition}

\begin{proof}
At $t=1$, the prolate equation \cref{eq:prolate-ode} gives
$\psi_0^{c\,\prime}(1)=\frac12(\lambda_0-c^2)\psi_0^c(1)$. Since
$\lambda_0\le c^2/3<c^2$ and $\psi_0^c(1)>0$ (\cref{sec:pswf}),
$\phi'(P/2^-)\ne0$. In the sum $S(y_\alpha/h)$, since $P$ is an integer,
the argument of one term enters the support at one end exactly when that
of another leaves it at the other, that is, when the particle crosses the
boundary of a stencil. The jumps of the two terms cancel there, but the derivative of
$S$ jumps by $2|\phi'(P/2^-)|$ (\cref{fig:charge}). Hence $S$ is not
constant, and $\sum_{m\ne0}|\widehat\phi(2\pi m)|^2>0$ by
\cref{eq:assigned-charge}. The quantity depends continuously on $c_w$,
which gives the positive minimum on a compact interval.
\end{proof}

Thus no choice of the bandwidth, the only parameter of the $\psi_0$
kernel, removes the variation of the assigned charge. The interval must
be compact, since, as $c_w\to0$, the $\psi_0$ kernel tends to $1/P$ times
the indicator function of $[-P/2,P/2]$, whose integer translates sum to
one.

With AD, charge conservation is not enough. For a kernel that conserves
charge, the charges that a unit charge at $\vy$ assigns to the grid have
a dipole moment whose components are $y_\alpha+hS_1(y_\alpha/h)$, where
\begin{equation}
 S_1\Bigl(\frac{y_\alpha}h\Bigr)=\sum_{\ell\in\Z}\Bigl(\ell-\frac{y_\alpha}h\Bigr)\,
   \phi\Bigl(\ell-\frac{y_\alpha}h\Bigr),
 \label{eq:assigned-dipole}
\end{equation}
whereas the particle has the dipole moment $\vy$. We say that the kernel
also conserves the dipole if $S_1(y_\alpha/h)=0$ for every $y_\alpha$, so
that, wherever the particle lies, the grid receives its exact dipole
moment.
The function $S_1$ has period one, and its $m$th Fourier coefficient is
$\int_{\R}t\phi(t)e^{2\pi\ii mt}\,dt=-\ii(d\widehat\phi/dq)(2\pi m)$, which
vanishes for $m=0$ since $\phi$ is even. For a kernel that conserves
charge but not the dipole, interpolating a linear mesh potential
$\mathbf E\cdot\vx$ with $\nabla\Phi$ gives a unit charge at $\vy$ the
force with components $-E_\alpha\bigl(1+S_1'(y_\alpha/h)\bigr)$ instead
of $-E_\alpha$. The spurious part $-E_\alpha S_1'(y_\alpha/h)$ does not
decrease with $h$, and the next proposition shows that it determines the
limit of the AD error.

\begin{proposition}
\label{prop:dipole-limit}
Let $\phi$ satisfy the hypotheses of \cref{prop:alias-limit}, be
continuous on $\R$, and have integer translates that sum to one. Then,
for every $\vk\in\Z^3\setminus\{\mathbf 0\}$,
\begin{equation}
 \begin{gathered}
 \lim_{h\to0}\eta_{\rm ik}(\vk)=0,\qquad
 \lim_{h\to0}\eta_{\rm AD}(\vk)=\frac D{1+D},\\
 D=\sum_{m\ne0}(2\pi m)^2\Bigl|\frac{d\widehat\phi}{dq}(2\pi m)\Bigr|^2 .
 \end{gathered}
 \label{eq:dipole-limit}
\end{equation}
If, in addition, $\phi$ conserves the dipole, then
$\eta_{\rm AD}(\vk)=O(h^2)$ as $h\to0$.
\end{proposition}

\begin{proof}
The functions $t\phi(t)$ and $t^2\phi(t)$ are continuous and piecewise
$C^2$, so integration by parts twice gives
$|(d\widehat\phi/dq)(q)|+|(d^2\widehat\phi/dq^2)(q)|\le C/(1+q^2)$
for a constant $C$. Since $\widehat\phi(2\pi m)=0$ for $m\ne0$, Taylor's
theorem gives
$\widehat\phi(q+2\pi m)=q\,(d\widehat\phi/dq)(2\pi m)+R_m(q)$ with
$|R_m(q)|\le Cq^2/(1+(2\pi|m|-\pi)^2)$ for $|q|\le\pi$. Hence
$\mathcal A(q)=O(q^2)$, so $T(\vk)=1+O(h^2)$ and
$\eta_{\rm ik}(\vk)=1-T(\vk)^{-2}=O(h^2)$. For AD, the bracket in
\cref{eq:T2-tensor} is
$h^{-2}\sum_{m\in\Z}(q+2\pi m)^2|\widehat\phi(q+2\pi m)/\widehat\phi(q)|^2$
with $q=h\xi_{\vk,\alpha}$. Its term with $m=0$ is $\xi_{\vk,\alpha}^2$.
The terms with $m\ne0$ are bounded by $C'\xi_{\vk,\alpha}^2/m^2$ for a
constant $C'$ and tend to
$\xi_{\vk,\alpha}^2(2\pi m)^2|(d\widehat\phi/dq)(2\pi m)|^2$, so, by
dominated convergence, the bracket tends to $\xi_{\vk,\alpha}^2(1+D)$.
Hence $T_2(\vk)\to|\vxi_\vk|^2(1+D)$ and
$\eta_{\rm AD}(\vk)\to1-1/(1+D)$. If $S_1\equiv0$, then
$(d\widehat\phi/dq)(2\pi m)=0$ for every $m$, only $R_m$ remains, and the
same bounds give the bracket $\xi_{\vk,\alpha}^2(1+O(h^2))$, so
$\eta_{\rm AD}(\vk)=O(h^2)$.
\end{proof}

By Parseval's identity, $D$ is the mean square over a grid cell of the
relative spurious force $S_1'(y_\alpha/h)$ in a linear potential. Thus with spectral differentiation, charge conservation suffices for the
aliasing error of each mode to vanish under refinement, while with AD,
a continuous kernel must also conserve the dipole. Both conditions together mean that the integer
translates of $\phi$ reproduce linear functions, which holds if and only
if $\widehat\phi$ vanishes to second order at every nonzero multiple of
$2\pi$~\cite{StrangFix}.

For the B-spline $B_P$, $\widehat\phi$ vanishes to order $P$ at every nonzero
multiple of $2\pi$, so both limits vanish, and
$\mathcal A(q)=\sum_{m\ne0}|q/(q+2\pi m)|^{2P}\sim2\zeta(2P)(q/2\pi)^{2P}$
as $q\to0$, where $\zeta$ is the Riemann zeta function. The relative
aliasing errors $\eta_{\rm ik}(\vk)^{1/2}$ and $\eta_{\rm AD}(\vk)^{1/2}$ of
$B_P$ at a fixed mode therefore decrease like $h^P$ and $h^{P-1}$,
respectively.

Since $\eta_{\rm ik}(\vk)$ and $\eta_{\rm AD}(\vk)$ depend on $\vk$ only
through $h\vxi_\vk$, these limits also describe the lowest modes of a
large box on a fixed grid, and a box of side $L$ has modes down to
$|\vxi_\vk|=2\pi/L$. The modes in a fixed interval of frequencies near
the bottom of the band carry a share of the force of the smooth part
that does not depend on $L$. For a kernel whose translates do not sum
to one, the proof of \cref{prop:alias-limit} shows that the relative
mean-square AD error of a mode grows like $1/(h|\vxi_\vk|)^2$ as $h|\vxi_\vk|$ decreases,
until it approaches one. Summed over the modes down to $2\pi/L$, these
errors make the relative mean-square AD error of the band grow linearly
with $L$ on a fixed grid, as long as the relative error of the lowest
modes stays well below one, which holds for the $\psi_0$ kernels of
\cref{sec:band-results} up to boxes of side $10^3r_c$ or more. For a kernel
that conserves charge and dipole, $\eta_{\rm AD}$ of the lowest modes
tends to zero instead (\cref{prop:dipole-limit}), and the error in a box
is close to its limit once the box holds many modes of the band.

\Cref{prop:alias-limit} applies to every SI kernel of fixed shape whose
integer translates do not sum to one, such as the truncated
Gaussian~\cite{SE,GSE} and Kaiser--Bessel~\cite{Gao} kernels. With AD, the aliasing
error of such a kernel at each mode tends to the force of the mode under
refinement at fixed order, and on a fixed grid, the error of such a
kernel grows with the box. For the AD error
to vanish under refinement and to depend little on the box, an SI kernel
of fixed shape must therefore assign to the grid the charge and, if the
kernel is continuous, also the dipole moment of the particle at every
position, as B-splines do. The SI
kernels of \cref{sec:si-design} do so exactly.

\subsection{The optimization problem}
\label{sec:si-design}

We turn now to the construction of the SI kernel. Molecular dynamics
codes keep the order $P$ fixed, often at four or five, and choose the
grid that meets the requested accuracy at the least cost~\cite{Ball,DH2}.
For an order $P$ and a tolerance $\varepsilon$, the best SI kernel is
therefore the one that meets $\varepsilon$ on the coarsest grid, that is,
with the smallest grid ratio $s$ in \cref{eq:grid-ratio}. The same
kernels give the lowest order that meets the tolerance on a given grid,
such as the grid that ESP uses, $s=1$, or an upsampled one: it is the lowest
order whose coarsest grid is no finer than the given grid, since, on finer
grids, the error of each optimized kernel decreases
(\cref{sec:band-results}).

We seek the best kernel among those with two further properties that
molecular dynamics needs, and we impose both by construction.
The first is charge and dipole conservation. By \cref{sec:refinement},
without it, the AD error of a continuous kernel of fixed shape does not vanish under
refinement, so that a tolerance below its smallest error calls for a
higher order; with it, every order meets every tolerance on a fine
enough grid. The second is smoothness: the computed energy must be twice
continuously differentiable where a particle crosses the boundary of a
stencil, as \cref{sec:split-smoothness} requires at the cutoff.

The kernels with both properties have a simple form. Let
$B_2(t)=\max(1-|t|,0)$ be the linear B-spline and $a=(P-2)/2$. For an
even function $g$ supported in $[-a,a]$, the kernel
\begin{equation}
 \phi=B_2*g,\qquad
 \widehat\phi(q)=\Bigl(\frac{\sin(q/2)}{q/2}\Bigr)^2\widehat g(q),
 \label{eq:si-factor}
\end{equation}
is supported in $[-P/2,P/2]$, and since $\phi''(t)=g(t+1)-2g(t)+g(t-1)$,
it is twice continuously differentiable on $\R$ if and only if $g$ is
continuous on $\R$, that is, $g(\pm a)=0$. Its transform vanishes to
second order at every nonzero multiple of $2\pi$, since that of $B_2$
does, so every kernel of this form conserves charge and dipole by
construction, whatever $g$ is, as the B-spline $B_P=B_2*B_{P-2}$ does.
The converse also holds.

\begin{proposition}
\label{prop:si-factor}
Every twice continuously differentiable SI kernel of order $P$ that
conserves charge and dipole has the form \cref{eq:si-factor}, with $g$
continuous and even.
\end{proposition}

\begin{proof}
The function $G(t)=\sum_{j\ge1}j\,\phi(t-j)$ is twice continuously
differentiable, satisfies $G(t+1)-2G(t)+G(t-1)=\phi(t)$, vanishes for
$t<-a$, and equals $tS(t)+S_1(t)=t$ for $t>a$; hence $g=G''$ is
continuous and supported in $[-a,a]$, and $\phi=B_2*g$. Since
$G(t)-G(-t)=tS(t)+S_1(t)=t$ for every $t$, $g$ is even.
\end{proof}

On the grid of spacing $h_0/s$ given by \cref{eq:grid-ratio}, we measure a kernel by the AD aliasing error of the
modes of the band $|\vxi_\vk|\le c_s/r_c$ relative to their force,
$(Q_{\rm AD}/\sum_\vk|\vxi_\vk|^2\Sh(\vxi_\vk)^2)^{1/2}$, with both sums
over these modes. This is the square root of the mean of
$\eta_{\rm AD}(\vk)$, weighted by the mean-square forces of the modes. Since
$\eta_{\rm AD}(\vk)$ depends on $\vk$ only through $h\vxi_\vk$ and is
bounded, the sums become integrals as $L\to\infty$ at fixed $h$, and the
limit
\begin{equation}
 E(\phi;s)=\Biggl(\frac{\int_{|\vxi|\le c_s/r_c}|\vxi|^2\Sh(\vxi)^2\,\eta_{\rm AD}(\vxi)\,d\vxi}
 {\int_{|\vxi|\le c_s/r_c}|\vxi|^2\Sh(\vxi)^2\,d\vxi}\Biggr)^{1/2},
 \label{eq:si-loss}
\end{equation}
in which $\eta_{\rm AD}(\vxi)$ is $\eta_{\rm AD}(\vk)$ with $\vxi_\vk$
replaced by $\vxi$, depends on the grid only through $s$ and not on the
box. We optimize the kernel for this limit, so that one kernel serves
every box; for a kernel that conserves charge and dipole, the limit is
close to the error in a box once the box holds many modes of the band
(\cref{sec:refinement}).

We call a kernel of order $P\ge4$ admissible if it has the form
\cref{eq:si-factor} with $g$ continuous and satisfies
\begin{equation}
 \begin{gathered}
 \int_{\R}\phi(t)\,dt=1,\qquad
 \phi'(t)\le0\quad(0<t<P/2),\\
 |\phi''(t)-\phi''(\tau)|\le|t-\tau|\max|B_P'''|\qquad(t,\tau\in\R),
 \end{gathered}
 \label{eq:si-admissible}
\end{equation}
where the maximum is over the intervals on which $B_P'''$ exists.

The conditions in \cref{eq:si-admissible} have the following roles. The
first normalizes the kernel. The sign condition makes $\phi$ a positive bump, as
$B_P$ and the $\psi_0$ kernel are, so that a particle assigns to every
grid point a charge of its own sign. The Lipschitz bound on $\phi''$
makes the smoothness of the kernel quantitative. The AD force is the
gradient of the interpolated potential, so its second derivatives along a
trajectory, which enter the local error of a Verlet step, are bounded in
terms of the Lipschitz constant of $\phi''$; we bound this constant by
that of $B_P$. The B-spline $B_P$ satisfies all the conditions.

The kernel that we seek is a solution of two nested problems. On the grid
of spacing $h_0/s$, the smallest error of an admissible kernel is
\begin{equation}
 E^*(s)=\min\bigl\{E(\phi;s):\ \phi\ \text{admissible}\bigr\},
 \label{eq:si-minimum}
\end{equation}
and for a tolerance $\varepsilon$, which fixes $c_s$ and the split, the
coarsest grid on which order $P$ meets the tolerance has the spacing
$h_0/s^*$, where
\begin{equation}
 s^*=\min\{s\ge1:\ E^*(s)\le\varepsilon\}.
 \label{eq:si-coarsest}
\end{equation}
We define the SI kernel of order $P$ for the tolerance $\varepsilon$ as
an admissible kernel that minimizes $E(\phi;s^*)$, call it the optimized
kernel of this order and tolerance, and denote it by $\phi_P^*$. In
words, among the admissible kernels of order $P$, which conserve charge
and dipole and are twice continuously differentiable, $\phi_P^*$ meets
the tolerance on the coarsest grid and is the most accurate on that grid.
A code uses the first grid of FFT-friendly size whose spacing
is at most $h_0/s^*$ (\cref{sec:parameters}), and on the finer grids of
\cref{sec:band-results}, up to $s=4$, the error of $\phi_P^*$ decreases. The
tolerance thus applies to the aliasing error of the band; the truncation
error of the split is at most the tolerance (\cref{sec:split-results}),
and \cref{sec:complete-results} measures the
complete error.

The Lipschitz bound on $\phi''$ is also what makes the minimum exist. The
aliasing error depends on $\phi$ only through the alias sums of
$\widehat\phi$ and $q\widehat\phi$, which are continuous in the norm of
$H^1(\R)$, whereas the values of $\phi'$ and $\phi''$ at $\pm P/2$ are
not. As with the split, kernels whose second derivative drops to zero in
ever thinner layers at the ends of the support can then approach an
infimum that no twice continuously differentiable kernel attains, and
their third derivatives grow without bound. The sign condition does not prevent this,
since these layers are monotone.

\begin{theorem}
\label{thm:si-existence}
Let $P\ge4$. For every $s\ge1$, the error $E(\phi;s)$ attains its
minimum $E^*(s)$ over the admissible kernels, and this minimum is at most
$E(B_P;s)$. For
every admissible kernel,
$\lim_{s\to\infty}E(\phi;s)=0$. Hence, for every $\varepsilon>0$, $s^*$ in
\cref{eq:si-coarsest} exists and is at most the smallest $s\ge1$ with
$E(B_P;s)\le\varepsilon$.
\end{theorem}

\begin{proof}
Let $M=\max|B_P'''|$. For an admissible kernel, $\phi''$ vanishes at
$\pm P/2$ and is $M$-Lipschitz, so $|\phi''(t)|\le M(P/2-|t|)$, and the
kernels and their first two derivatives are uniformly bounded and
equicontinuous. By the Arzel\`a--Ascoli theorem, every sequence of
admissible kernels has a subsequence that converges uniformly on $\R$
together with its first two derivatives. The limit satisfies
\cref{eq:si-admissible}, and its transform keeps the double zeros at the
nonzero multiples of $2\pi$, so, by \cref{prop:si-factor}, it has the form
\cref{eq:si-factor}: the admissible kernels form a compact set in
$C^2(\R)$. Three integrations by parts give $|\widehat\phi(q)|\le MP/|q|^3$,
so the alias sums $\sum_\vm\widehat\Phi(\vxi_\vm)^2$ and
$\sum_\vm|\vxi_\vm|^2\widehat\Phi(\vxi_\vm)^2$ converge uniformly over the
admissible kernels and depend continuously on the kernel and on
$h\vxi$. By \cref{eq:alias-ratios,eq:Q-AD-ik}, $\eta_{\rm AD}(\vxi)$ is one
minus $|\vxi|^2\widehat\Phi(\vxi)^4$ divided by the product of these sums,
a number in $[0,1]$ that depends continuously on the kernel wherever both
sums are positive. For a given kernel, they vanish only on finitely many
planes $\xi_\alpha=\text{const}$ in the band, since
$\sum_{m}\widehat\phi(q+2\pi m)^2$ and
$\sum_{m}(q+2\pi m)^2\widehat\phi(q+2\pi m)^2$ are real-analytic functions
of $q$ ($\widehat\phi$ is entire and bounded by $MPe^{P|\Im z|/2}/|z|^3$),
whose integrals over a period are $2\pi\int_{\R}\phi(t)^2\,dt$ and
$2\pi\int_{\R}\phi'(t)^2\,dt$, respectively, both positive. By dominated convergence,
$E(\phi;s)$ is continuous on this compact set and attains its minimum,
and $B_P$ is admissible, which gives $E^*(s)\le E(B_P;s)$. For a fixed
admissible kernel, \cref{prop:dipole-limit} gives $\eta_{\rm AD}\to0$ at
every point of the band as $s\to\infty$, and dominated convergence gives
$E(\phi;s)\to0$. Finally, let $s_n\to s$ with $E^*(s_n)\le\varepsilon$,
and let $\phi_n$ attain $E^*(s_n)$. A subsequence converges in $C^2(\R)$
to an admissible kernel $\phi$, and $E(\phi_n;s_n)\to E(\phi;s)$ by the
same argument, so $E^*(s)\le\varepsilon$. The set in
\cref{eq:si-coarsest} is therefore closed, and it is not empty since
$E^*(s)\le E(B_P;s)\to0$.
\end{proof}

In other words, the optimized kernel of order $P$ never needs a finer grid than the
B-spline of the same order, and its error tends to zero under refinement,
whereas that of the $\psi_0$ kernel does not. Since the form \cref{eq:si-factor}
imposes only the double zeros of $\widehat\phi$ at the nonzero multiples
of $2\pi$, this decrease is in general only of first order in $h$
(\cref{prop:dipole-limit}), compared with order $P-1$ for $B_P$. The
optimized kernel thus trades the rate for accuracy on its coarsest grid,
and a tighter tolerance calls for the kernel designed for it, not only
for a finer grid.

\subsection{Construction of the kernels}
\label{sec:si-construction}

In order to solve the problems \cref{eq:si-minimum,eq:si-coarsest}
numerically, we expand $g$ in even PSWFs on $[-a,a]$,
\begin{equation}
 g(t)=\sum_{j=0}^{K-1}b_{2j}\,\psi_{2j}^{c_w}(t/a),\qquad |t|\le a,
 \label{eq:si-expansion}
\end{equation}
where the bandwidth $c_w$ only selects the basis; we take that of the
$\psi_0$ kernel of \cref{sec:band-results} on the grid with $s=1$. The equality
constraints, $\int_{\R}\phi(t)\,dt=1$ and $g(a)=0$, are linear in the
coefficients, and we impose the inequalities of \cref{eq:si-admissible}
at Chebyshev points of the intervals on which $\phi'''$ is smooth. The
computation then has two steps.
\begin{enumerate}[leftmargin=*]
\item \emph{The minimum on one grid.} On the grid of spacing $h_0/s$, we minimize
$E(\phi;s)$ over the coefficients by sequential quadratic
programming~\cite{NocedalWright}, with the integral \cref{eq:si-loss}
evaluated on a lattice in the band. By
\cref{eq:alias-power,eq:T2-tensor}, the error depends on the kernel only
through the one-dimensional alias sums at the coordinates of the lattice,
fewer than a hundred values of $q$. Since $\widehat\phi$ is linear in the
coefficients, the cell integrals \cref{eq:alias-cell} of the basis
functions are computed once for each grid, and an evaluation of the
error and its gradient takes a few hundredths of a second. The problem
is not convex, and the method finds local minima. We start it from $B_P$
and from the kernel with $g$ a multiple of $\psi_0^{c_w}(t/a)$, both
expressed in the expansion, and from the minimizer on the previous grid,
and keep the smallest minimum.
\item \emph{The coarsest grid.} We compute the minimum with $K=64$ on the
grids with $s=1.1^j$, $j=0,1,\ldots$, until it is at most $\varepsilon$, and
then bisect the last interval until the computed $s^*$ is known to within
$0.1\%$. The computed minima decrease with $s$ for every order and
tolerance, so this search finds the smallest ratio at which the computed
minimum meets the tolerance.
\end{enumerate}
Each minimization takes a few seconds on one core, and the search for
$s^*$ about a minute for each order and tolerance; the kernels are
computed once and tabulated (\cref{sec:parameters}). In the rest of the
paper, $\phi_P^*$ denotes the computed kernel.

\subsection{Choosing the parameters}
\label{sec:parameters}

For the Coulomb interaction, a requested tolerance $\varepsilon$, a cutoff
$r_c$, and an order $P$ that the code supports determine the parameters of
the method. The tolerance fixes the splitting bandwidth $c_s$: the
bandwidth of ESP's rule $\chi(1)=\varepsilon$, raised at
$10^{-2}$ and $10^{-3}$ so that the truncation error of the $C^2$ split
meets the tolerance (\cref{sec:split-results}). The formula
$c_s=\ln(1/\varepsilon)+0.86\ln\ln(1/\varepsilon)+0.91$ reproduces the
bandwidth of ESP's rule to within $0.01$ for the tolerances from
$10^{-2}$ to $10^{-8}$. The tolerance and the order fix the kernel $\phi_P^*$ and
its coarsest ratio $s^*$ \cref{eq:si-coarsest}, and hence the grid
spacing $h=\pi r_c/(c_ss^*)$; \cref{tab:coarsest} gives $c_s$ and $h/r_c$
for $P=4$, $5$, and $6$. For a box of side $L$, a code takes the first
grid that has at least $L/h$ points in each dimension and a size that the
FFT handles efficiently, and the method uses $\phi_P^*$ with its optimal influence function
\cref{eq:I-opt} on this grid. A code that
uses the coarsest grid of $B_P$, which is never coarser than that of
$\phi_P^*$ (\cref{tab:coarsest}), can keep this grid with $\phi_P^*$,
whose error decreases on finer grids (\cref{sec:band-results}). With these
parameters, the complete errors of \cref{sec:complete-results}, at the
tolerance $10^{-4}$, stay below the tolerance, relative to the RMS force,
for independent particles and for liquid water.

The choice among orders depends on the hardware: it balances the cost of
spreading and interpolation, which grows like $P^3$ per particle, against
the cost of the FFTs, which grows with the number of grid points. A simple rule gives one order for
every tolerance from $10^{-2}$ to $10^{-8}$: in the computations of
\cref{sec:band-results}, $\phi_P^*$ with $P=\log_{10}(1/\varepsilon)+2$
meets the tolerance on the grid of spacing $\pi r_c/(1.1c_s)$, on which
neither $B_P$ nor the $\psi_0$ kernel of the same order does.

The kernel $\phi_P^*$ has no closed form, so a code evaluates it from a
precomputed table, which is how ESP evaluates the $\psi_0$
kernel~\cite{ESP}. We divide each
grid interval of the support into $128$ equal parts and tabulate
$\phi_P^*$ and its first two derivatives at the ends of the parts. On
each part, $\phi_P^*$ is evaluated by the polynomial of degree five that
matches these values at both ends, and its derivative, which AD needs,
by the derivative of this polynomial. This table is accurate enough: for
every order and tolerance of \cref{sec:band-results}, it changes
$E(\phi_P^*;s^*)$ by less than $0.4\%$. These tables, the coefficients of
the kernels, and programs that evaluate them are
available~\cite{PMEkernels}.

\FloatBarrier
\section{Numerical results}
\label{sec:numerics}

We test the $C^2$ splits of \cref{sec:split-design} for $\nu=1,\ldots,6$,
then the SI kernels of \cref{sec:si} on the band of the $C^2$ Coulomb
split and in complete Coulomb calculations, and finally the conservation
of energy with the two kernels together over $10^8$ time steps. The
programs and data of \cref{sec:split-results,sec:band-results} are
available~\cite{PMEkernels}.

\subsection{The \texorpdfstring{$C^2$}{C2} splitting kernels}
\label{sec:split-results}

We computed the $C^2$ splits of \cref{sec:split-construction} for the
tolerances $\varepsilon=10^{-2},10^{-3},\ldots,10^{-8}$, with one
bandwidth for all $\nu$, as in~\cite{MLIP}: the bandwidth of ESP's rule
$\chi(1)=\varepsilon$, raised by $6\%$ at $10^{-2}$ and by $0.3\%$ at
$10^{-3}$ so that the truncation error of the $C^2$ Coulomb split meets
the tolerance. We keep calling the grid with $s=1$ the grid that ESP
uses, although at these two tolerances it is finer than ESP's own grid
by the same factors.
For the Coulomb interaction, the relative truncation error of the $C^2$
split is then at most the tolerance, and for the higher powers, it lies
between about one and eight times the tolerance. It grows with the power,
as the jumps of the $\psi_0$ splits at the cutoff do
(\cref{sec:split-smoothness}), so that with either split, the higher
powers need a larger bandwidth to meet the tolerance. To reach the
truncation error of the $\psi_0$ split in a box of side $10r_c$, with all
the modes of the grid retained, the $C^2$ Coulomb split needs grids with
$9\%$ to $14\%$ more points in each dimension than the grid that ESP
uses with its own bandwidth. This mild refinement is the price of
smoothness at the cutoff.

\subsection{The SI kernels}
\label{sec:band-results}

We computed the optimized kernels $\phi_P^*$ of
\cref{sec:si-construction} for the tolerances
$\varepsilon=10^{-2},\ldots,10^{-8}$ and the orders $P=4,\ldots,10$, with
the $C^2$ Coulomb split and $1\le s\le4$, and compared them with $B_P$
and with the $\psi_0$ kernel \cref{eq:psi0-kernel}, whose bandwidth $c_w$
we chose on each grid to minimize its error. Every kernel uses its own
optimal influence function \cref{eq:I-opt}, so the B-spline baseline is
the P3M method with AD~\cite{Ball}. The errors of $B_P$ and $\phi_P^*$ are
those of \cref{eq:si-loss}, in the limit of large boxes; since both
kernels conserve charge and dipole, their errors in a box of side $10r_c$
differ little from this limit (\cref{sec:refinement}). The error of the
$\psi_0$ kernel grows with the box, and we give it in a box of side
$10r_c$, comparable to the smaller systems in the benchmarks
of~\cite{ESP}, with $c_w$ chosen for this box. We write $s^*$
also for the coarsest ratios of $B_P$ and of the $\psi_0$ kernel.

\begin{figure}[htbp]
\centering
\includegraphics[width=0.93\linewidth]{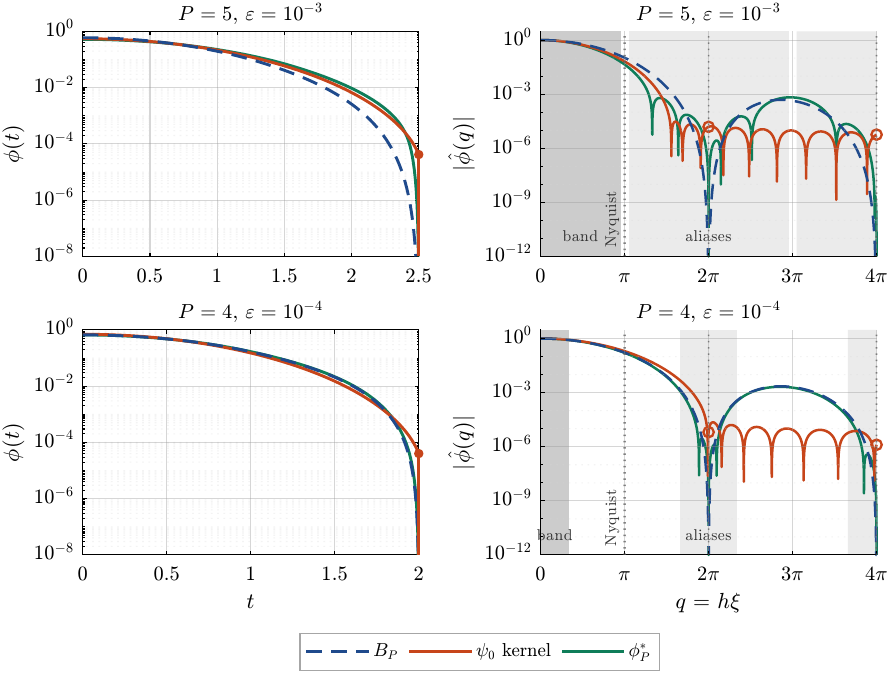}
\caption{The SI kernels $\phi(t)$ (left) and their transforms
$|\widehat\phi(q)|$, $q=h\xi$ (right), for $P=5$ at the tolerance
$10^{-3}$ (top) and $P=4$ at $10^{-4}$ (bottom), on the coarsest grid of
$\phi_P^*$: $B_P$, the $\psi_0$ kernel for this grid, and $\phi_P^*$, each
with unit integral. In the transforms, the band $|q|\le\pi/s$ is shaded
dark and its aliases around $2\pi$ and $4\pi$ light, and the dotted lines
mark the Nyquist frequency and $2\pi$ and $4\pi$. By construction,
$\phi_P^*$ shares two properties with $B_P$: it is twice continuously
differentiable, and it conserves charge; that is, its transform vanishes
at the nonzero multiples of $2\pi$. The $\psi_0$ kernel jumps to zero at
the end of its support (dot), and its transform does not vanish at
$2\pi$ and $4\pi$ (circles).}
\label{fig:kernels}
\end{figure}

\Cref{fig:kernels} shows the three kernels for the default orders of
molecular dynamics codes, $P=5$ at $10^{-3}$ and $P=4$ at $10^{-4}$, on
the coarsest grid of $\phi_P^*$. Their transforms show the charge
conservation that \cref{sec:refinement} requires and \cref{sec:si-design}
imposes on $\phi_P^*$: those of $\phi_P^*$ and $B_P$ vanish at the
nonzero multiples of $2\pi$, where that of the $\psi_0$ kernel does not. At the alias $2\pi-\pi/s$ of the band edge, the
transform of $\phi_P^*$ is the smallest of the three.

On the grid that ESP uses, the coarsest of the study, $\phi_P^*$ needs
one order less than the $\psi_0$ kernel at $10^{-3}$ and $10^{-4}$, two
less at $10^{-5}$, and the same order at $10^{-2}$
(\cref{tab:minimal-order}); at $10^{-6}$, it meets the tolerance with
$P=9$, whereas neither of the other kernels does with an order up to ten.
B-splines need one and four orders more than $\phi_P^*$ at $10^{-2}$ and
$10^{-3}$, respectively, and an order above ten from $10^{-4}$ on. Below $10^{-6}$, no
kernel of order up to ten meets the tolerance on this grid, whose band
reaches the Nyquist frequency along the axes, where no kernel can tell a
mode from its alias; $\phi_P^*$ meets $10^{-7}$ with $P=10$ on a grid
$0.8\%$ finer.

\begin{table}[htbp]
\centering
\caption{The smallest order $P\le10$ whose SI kernel meets the tolerance
$\varepsilon$ on the grid that ESP uses, $s=1$, in the relative AD
aliasing error with the $C^2$ Coulomb split: $B_P$ and $\phi_P^*$ in the
limit of large boxes \cref{eq:si-loss}, and the $\psi_0$ kernel in a box
of side $10r_c$. Below $10^{-6}$, no kernel of order up to ten meets the
tolerance on this grid.}
\label{tab:minimal-order}
\small
\begin{tabular}{@{}lccc@{}}
\toprule
$\varepsilon$ & $B_P$ & $\psi_0$ & $\phi_P^*$\\
\midrule
$10^{-2}$ & 6 & 5 & 5\\
$10^{-3}$ & 10 & 7 & 6\\
$10^{-4}$ & $>10$ & 8 & 7\\
$10^{-5}$ & $>10$ & 10 & 8\\
$10^{-6}$ & $>10$ & $>10$ & 9\\
\bottomrule
\end{tabular}
\end{table}

Molecular dynamics codes usually keep the order at four or five and
refine the grid instead~\cite{GromacsGPU,ESP}. \Cref{tab:coarsest} gives,
for $P=4$, $5$, and $6$, the largest grid spacing with which each kernel
meets the tolerance, the quantity that a code sets. In all 49 cases of
the study, no kernel meets a tolerance on a coarser grid than $\phi_P^*$;
for $B_P$, \cref{thm:si-existence} guarantees this for the exact
minimizers. Since the number of grid points is proportional to $h^{-3}$,
$B_P$ needs up to $4.8$ times as many grid points as $\phi_P^*$, and up
to $3.5$ times as many with the default orders four and five. With these
orders and $s\le4$, the $\psi_0$ kernel meets no tolerance below
$10^{-3}$, whereas $\phi_P^*$ meets $10^{-4}$ with $P=4$ and $10^{-6}$
with $P=5$; over all orders, the $\psi_0$ kernel needs up to $14$ times
as many grid points as $\phi_P^*$.

\begin{table}[htbp]
\centering
\caption{The largest grid spacing $h/r_c=\pi/(c_ss^*)$ with which the SI
kernel of order $P$ meets the tolerance $\varepsilon$, with the $C^2$
Coulomb split of bandwidth $c_s$: $B_P$ and $\phi_P^*$ in the limit of
large boxes, and the $\psi_0$ kernel in a box of side $10r_c$. A dash
means that no grid with $s\le4$ suffices.}
\label{tab:coarsest}
\footnotesize
\setlength{\tabcolsep}{3pt}
\begin{tabular}{@{}lcccccccccc@{}}
\toprule
 & & \multicolumn{3}{c}{$P=4$} & \multicolumn{3}{c}{$P=5$} & \multicolumn{3}{c}{$P=6$}\\
\cmidrule(lr){3-5}\cmidrule(lr){6-8}\cmidrule(lr){9-11}
$\varepsilon$ & $c_s$ & $B_P$ & $\psi_0$ & $\phi_P^*$ & $B_P$ & $\psi_0$ & $\phi_P^*$ & $B_P$ & $\psi_0$ & $\phi_P^*$\\
\midrule
$10^{-2}$ & $7.22$ & $0.359$ & $0.365$ & $0.407$ & $0.412$ & $0.435$ & $0.435$ & $0.435$ & $0.435$ & $0.435$\\
$10^{-3}$ & $9.50$ & $0.175$ & -- & $0.205$ & $0.246$ & $0.263$ & $0.315$ & $0.285$ & $0.327$ & $0.331$\\
$10^{-4}$ & $12.02$ & $0.0771$ & -- & $0.0894$ & $0.139$ & -- & $0.199$ & $0.183$ & $0.2$ & $0.248$\\
$10^{-5}$ & $14.52$ & -- & -- & -- & $0.0754$ & -- & $0.115$ & $0.115$ & -- & $0.171$\\
$10^{-6}$ & $16.98$ & -- & -- & -- & -- & -- & $0.0589$ & $0.0704$ & -- & $0.112$\\
$10^{-7}$ & $19.42$ & -- & -- & -- & -- & -- & -- & $0.0427$ & -- & $0.072$\\
$10^{-8}$ & $21.84$ & -- & -- & -- & -- & -- & -- & -- & -- & $0.0451$\\
\bottomrule
\end{tabular}
\end{table}

\begin{figure}[htbp]
\centering
\includegraphics[width=0.95\linewidth]{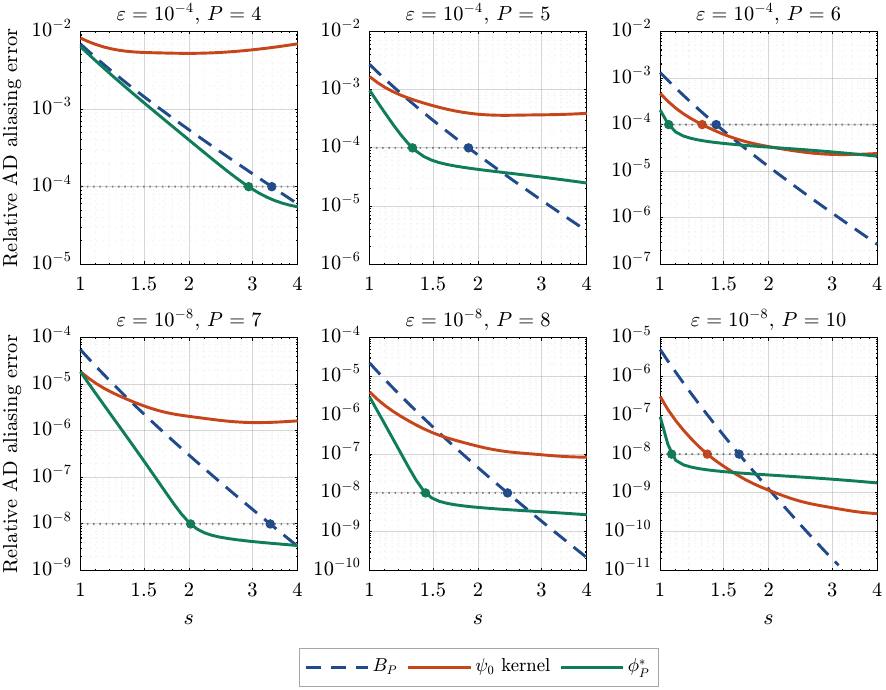}
\caption{Relative AD aliasing error against the grid ratio $s$, for the
tolerances $10^{-4}$ (top) and $10^{-8}$ (bottom), with the $C^2$ Coulomb
split: $B_P$ and $\phi_P^*$ in the limit of large boxes
\cref{eq:si-loss}, and the $\psi_0$ kernel in a box of side $10r_c$,
whose error grows in larger boxes. Errors below $10^{-11}$, near the
floor of their computation, are not shown. The dotted lines mark the
tolerance, and the dots the coarsest ratios $s^*$. In every panel,
$\phi_P^*$ meets the tolerance on a coarser grid than the other two
kernels, and its error keeps decreasing on finer grids.}
\label{fig:coarsest}
\end{figure}

\Cref{fig:coarsest} shows the errors under refinement, for a typical
tolerance and the tightest one of the study. The error of $\phi_P^*$
tends to zero under refinement, as \cref{prop:dipole-limit} predicts,
and for every order and tolerance of the study, it decreases from $s^*$
to $s=4$. It decreases more slowly than that of $B_P$, since $\phi_P^*$
trades the rate for accuracy on its coarsest grid (\cref{sec:si-design}),
but the error of $B_P$ falls below that of $\phi_P^*$ only on grids finer
than the coarsest grid of $B_P$. At $10^{-4}$,
the error of the $\psi_0$ kernel levels off above the tolerance for $P=4$
and $5$. On a fixed grid, its error grows with the box, as
\cref{sec:refinement} predicts: on the coarsest grids of $\phi_P^*$ for $P=4$ and $6$ at
$10^{-4}$ and for $P=6$ at $10^{-6}$, it grows by factors of $1.5$ to $2.4$
from $L=10r_c$ to $L=80r_c$. With
$P=10$ at $10^{-8}$, in contrast, the error of the $\psi_0$ kernel falls
below that of $\phi_P^*$ on fine grids, since its bandwidth adapts to
each grid while $\phi_P^*$ is designed for its coarsest grid.

\subsection{Complete Coulomb calculations}
\label{sec:complete-results}

The comparison of \cref{sec:band-results} excludes the truncation error
and the self force. We now use the kernels designed for the tolerance
$10^{-4}$, typical in molecular dynamics~\cite{ESP}, unchanged in
complete AD calculations of the periodic Coulomb forces. The calculations
include the exact local part, retain all grid frequencies with
$|k_\alpha|<n_f/2$, and use the optimal influence function
\cref{eq:I-opt} for the full reference. We compute the forces of the
three kernels on the coarsest grid of each kernel (the first grid of the
box with even $n_f$ and $s\ge s^*$), for each order from four, the
default of GROMACS on GPUs, to seven, the lowest with which $\phi_P^*$
meets the tolerance on the grid that ESP uses. Against a Gaussian Ewald
sum with a relative accuracy of about $10^{-15}$, we report the charge-scaled RMS
force error
\begin{equation}
 \mathcal E_{\rm pair}=\Biggl[\frac1J\sum_{j=1}^J
   \frac{\sum_i|\mathbf F_i^{\rm mesh}+q_i\nabla u^s(\vx_i)|^2}{M_2^2-M_4}\Biggr]^{1/2}
 \label{eq:sample-pair-rms}
\end{equation}
over $J$ samples, where the inner sum and the charge moments refer to
sample $j$; the same expression with the exact forces in place of the
errors, the RMS force, sets the scale.

For independent particles, $J=16$ configurations of 64 uniformly
distributed unit charges of both signs in a box of side $10r_c$ confirm
the error model. With each particle's self force subtracted,
$\mathcal E_{\rm pair}^2$ estimates $\sigma_{\rm pair,F}^2$
(\cref{prop:force-variance}), and on its own coarsest grid, each kernel
has an error of about $7\times10^{-5}$ of the RMS force or less, as
\cref{eq:si-loss} predicts. On the coarsest grid of $\phi_P^*$, the errors of $B_P$
and of the $\psi_0$ kernel are $1.4$ to $8.2$ and $1.9$ to $52$ times
that of $\phi_P^*$, respectively, and the paired-bootstrap $95\%$
confidence intervals of these ratios exclude one. The kernel $\phi_P^*$ has the
smallest error on every other grid as well, except on the coarsest grid
of $B_7$, where the error of the $\psi_0$ kernel is $3\%$ smaller, within
the statistical uncertainty. It also has the smallest self force, by
factors of $3.1$ to $240$. By \cref{eq:force-variance}, the self term of
$N$ unit charges is the smaller one when
$N-1>\sigma_{\rm self,F}^2/\sigma_{\rm pair,F}^2$, and since
$\sigma_{\rm pair,F}^2$ falls like the inverse volume of the box at a
fixed grid spacing, this holds for $\phi_P^*$ from about $0.003$ charges
per grid point on.

The test that matters in practice is a molecular liquid, whose charges
are correlated, so that \cref{prop:force-variance} does not apply. We
take an equilibrated configuration of 901 SPC/E water
molecules~\cite{SPCE} from a molecular dynamics simulation, with 2,703
charges in a cubic box of side 30~\AA, and use $r_c=10$~\AA{} as the unit
of length, so that $L=3r_c$. For each kernel, the grid is the first grid
of this box with even $n_f$ whose ratio reaches the coarsest ratio of the
kernel; these grids have $0.04$ to $1.6$ charges per grid point. We keep
the self forces, as molecular dynamics codes do, and average over $J=16$
random translations, which sample the positions of the molecules
relative to the grid. On the coarsest grid of $\phi_P^*$
(\cref{tab:water}), the error of $\phi_P^*$ is at most $9\times10^{-6}$
of the RMS force, and those of $B_P$ and of the $\psi_0$ kernel are $1.7$
to $6.6$ and $1.03$ to $14$ times as large, respectively; with $P=7$,
$\phi_P^*$ comes within about $8\%$ of the truncation error of the grid.
It is also the most accurate kernel on the coarsest grids of $B_P$ and of
the $\psi_0$ kernel, except on that of $B_7$, where the error of the
$\psi_0$ kernel is $3\%$ smaller. Subtracting the self forces would
increase the errors in water, by factors of up to $3.3$ for the $\psi_0$
kernel and $1.5$ for $\phi_P^*$, since the self forces largely cancel the
mesh errors of the interactions within each molecule.

\begin{table}[htbp]
\centering
\caption{Complete Coulomb AD force errors \cref{eq:sample-pair-rms}, with
the self forces kept, for an equilibrated configuration of 901 SPC/E water
molecules in a box of side $3r_c$, averaged over 16 translations, with the
$C^2$ Coulomb split and the kernels for the tolerance $10^{-4}$, on the
coarsest grid of $\phi_P^*$, of spacing $h$, and the truncation error of
the grid. The RMS force is $2.83$.}
\label{tab:water}
\small
\begin{tabular}{@{}rccccc@{}}
\toprule
$P$ & $h/r_c$ & $B_P$ & $\psi_0$ & $\phi_P^*$ & Truncation\\
\midrule
4 & $0.0882$ & $2.64\times10^{-5}$ & $2.12\times10^{-4}$ & $1.51\times10^{-5}$ & $2.91\times10^{-6}$\\
5 & $0.188$ & $1.65\times10^{-4}$ & $8.37\times10^{-5}$ & $2.49\times10^{-5}$ & $1.36\times10^{-5}$\\
6 & $0.214$ & $9.08\times10^{-5}$ & $2.26\times10^{-5}$ & $1.84\times10^{-5}$ & $1.50\times10^{-5}$\\
7 & $0.25$ & $8.37\times10^{-5}$ & $1.88\times10^{-5}$ & $1.81\times10^{-5}$ & $1.68\times10^{-5}$\\
\bottomrule
\end{tabular}
\end{table}

\subsection{Long-time energy conservation}
\label{sec:nve}

\begin{figure}[htbp]
\centering
\includegraphics[width=\linewidth]{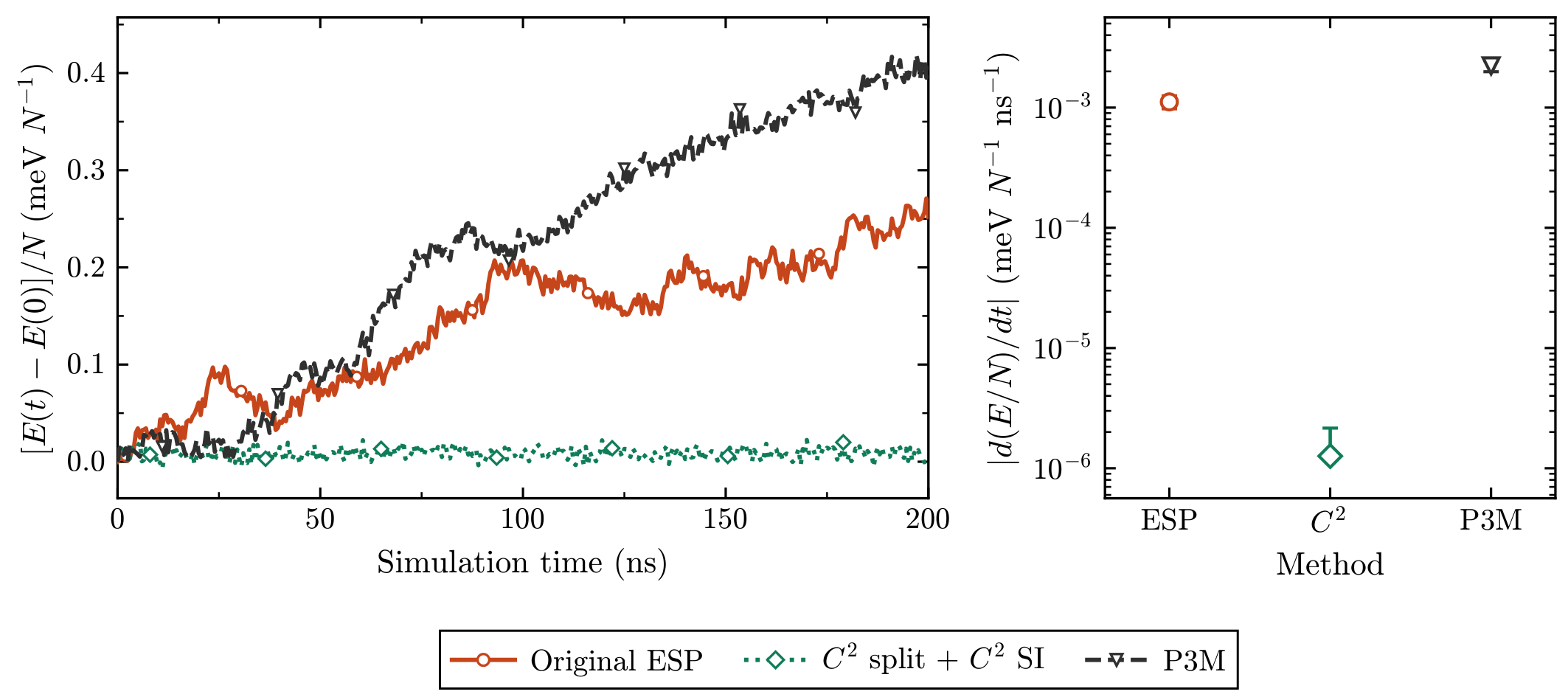}
\caption{Energy conservation in simulations of water without a
thermostat, with ESP, P3M, and the $C^2$ scheme, whose splitting
and SI kernels are twice continuously differentiable. Left: the deviation
of the total energy per atom from its initial value over $200$~ns.
Right: the magnitudes of the drift rates, estimated by least-squares fits
of linear functions of time to the curves on the left, with approximate
$95\%$ confidence intervals; the lower end of the interval for the $C^2$
scheme lies below the plotted range. Only with the $C^2$ scheme does the
energy stay close to its initial value.}
\label{fig:nve}
\end{figure}

To see how the smoothness of the kernels affects the conservation of
energy over long times, we
simulated $512$ rigid SPC/E water molecules~\cite{SPCE} in a periodic box
of side $24.83$~\AA, at a density of $1.0$~g\,cm$^{-3}$, without a
thermostat, so that the exact dynamics conserves the total energy and
any drift of the computed energy is a numerical error.
All simulations started from the same positions and velocities, at an
initial temperature of $300$~K, and used the velocity Verlet method with
a time step of $2$~fs, with the rigid geometry of the molecules
maintained by RATTLE~\cite{RATTLE} to a tolerance of $10^{-12}$ in at
most $500$ iterations. We used a modified version of LAMMPS (4~July
2026)~\cite{LAMMPS}. The Lennard-Jones interaction was switched off
between $7$ and $9$~\AA\ by a $C^3$ switching function. We compare
ESP~\cite{ESP}, P3M~\cite{HE}, and the $C^2$ scheme, which uses the
$C^2$ Coulomb split and the optimized SI kernel $\phi_5^*$ for the
tolerance $10^{-4}$. All three use the cutoff $9$~\AA, the order $P=5$,
and AD. ESP and the $C^2$ scheme use $15^3$ grids and $c_s=12.021$, with
$c_w=12.91$ and $12.43$, respectively. P3M uses the B-spline $B_5$ and
the parameters that LAMMPS chooses for the tolerance $10^{-4}$: the
Gaussian split with the Ewald parameter $\alpha=0.33$~\AA$^{-1}$ and a
$30^3$ grid, which has eight times as many points as the $15^3$ grid.
Measured against a
converged Ewald sum on the same $200$ configurations, sampled at
intervals of $1$~ns, the relative RMS force errors are
$1.12\times10^{-4}$ for ESP, $1.03\times10^{-4}$ for P3M, and
$4.08\times10^{-5}$ for the $C^2$ scheme, the most accurate of the
three.

\Cref{fig:nve} shows the deviation of the energy, recorded every
$20$~fs over $200$~ns, that is, $10^8$ time steps, and the magnitudes of
the drift rates, obtained by fitting a linear function of time to the signed
deviations. Since the deviations are correlated in time, we sum the
products of the residuals with the centered times over $20$ consecutive
blocks of about $10$~ns and estimate the variance of each slope from the
block sums by the Newey--West estimator~\cite{NeweyWest}, with Bartlett
weights up to lag $2$ and a finite-sample correction; the error bars are
the resulting approximate $95\%$ confidence intervals. The energy drifts
steadily with ESP and P3M, whereas it stays close to its initial value
with the $C^2$ scheme. The fitted drift rates are $1.12\times10^{-3}$,
$2.22\times10^{-3}$, and $1.26\times10^{-6}$~meV\,$N^{-1}$\,ns$^{-1}$
for ESP, P3M, and the $C^2$ scheme, respectively. The rate of the $C^2$
scheme is nearly three orders of magnitude smaller than that of ESP,
although the force errors of the three schemes are of the same order.
This rate is comparable to the uncertainty of its estimate. With the
heat capacity of liquid water, about $75$~J\,mol$^{-1}$\,K$^{-1}$, the
energy gained in $200$~ns would heat the water by about $0.9$~K with ESP
and $1.7$~K with P3M, and by about $1$~mK with the $C^2$ scheme. The
three simulations differ only in the treatment of the electrostatic
interactions, with the same constraints and the same switching of the
Lennard-Jones interaction, so the excess drift with ESP and P3M comes
from their kernels.

\FloatBarrier
\section{Conclusions}
\label{sec:conclusions}

We have presented a framework for designing the two kernels of PME
methods for the needs of molecular dynamics. The smoothness that long
molecular dynamics simulations require is imposed on both kernels, and on the SI kernel
also the conservation of charge and dipole that upsampling with analytic
differentiation requires; within these constraints, each kernel is
chosen to make the error of the method small. The resulting kernels
improve on both the Gaussian split with B-splines, used in PME methods
for three decades, and the recent PSWF kernels of ESP.

Our splits give every inverse power that we treat, including the Coulomb
and the dispersion interaction, a local part that is smooth at the
cutoff to the order that the time integration needs. Since their window
functions are combinations of the most concentrated even PSWFs that admit
these smoothness conditions, they stay close to the optimal Fourier
concentration, and the grid becomes only mildly finer: for the
Coulomb interaction, our twice continuously differentiable split reaches
the truncation error of the $\psi_0$ split of ESP with at most $14\%$
more grid points along each axis.

For the SI kernel, molecular dynamics needs a kernel of low order that
is accurate on coarse grids and whose error falls as the grid is
refined. With
analytic differentiation, which most codes use, refining the grid makes
the aliasing error of a kernel of fixed shape vanish only if the kernel
conserves charge and, if it is continuous, also the dipole. The $\psi_0$
kernel and the truncated Gaussian and Kaiser--Bessel kernels do not
conserve charge, which explains why they gain little from upsampling. Our SI kernels conserve both by construction, are twice
continuously differentiable, and are optimized on the coarsest grid on
which their order meets the tolerance. They never need a finer grid than
the B-spline of equal order, and with the default orders four and five,
they cut the number of grid points of B-splines by factors of up to
$3.5$. On the coarse grids of ESP, their order is lower than that of
the $\psi_0$ kernel at most tolerances, and never higher. The gains carry over to complete calculations for liquid water.

Together, the new kernels make the energy a smooth function of the
positions of the particles, and when liquid water is simulated without
a thermostat, they slow the drift of the energy by about three orders
of magnitude relative to the kernels of P3M and ESP. The SI kernels
depend only on the order and the tolerance, not on the grid or the box,
so they are tabulated once; existing codes adopt the new kernels
without changing the structure of their PME methods. Since PME methods
are the default for electrostatics in the major molecular dynamics
packages, the new kernels can make many simulations both cheaper,
through smaller grids, and more reliable, through better conservation of
energy.

\section*{Acknowledgments}
We thank Yutong Zhao at NVIDIA, and Alex Barnett and Leslie Greengard at
the Flatiron Institute, for useful discussions. The work of J. Liang was
supported by the National Natural Science Foundation of China (grant
12401570) and the China Postdoctoral Science Foundation (grant
2024M751948).

\bibliographystyle{siamplain}
\bibliography{pme_kernels}

\end{document}